\documentclass[a4paper,12pt]{amsart}
\usepackage{amssymb,amscd}
\usepackage{calrsfs}

\theoremstyle{plain}
\newtheorem{theo}{Theorem}[section]
\newtheorem{prop}[theo]{Proposition}
\newtheorem{lemm}[theo]{Lemma}
\newtheorem{coro}[theo]{Corollary}

\theoremstyle{definition}

\newtheorem*{exam}{Example}

\theoremstyle{remark}
\newtheorem{rema}[theo]{Remark}
\newtheorem*{note}{Note}

\renewcommand{\u}{u_i^I}

\newcommand{\V}{\mathcal V}
\newcommand{\E}{\mathcal E}

\newcommand{\field}[1]{\mathbb{#1}}
\newcommand{\C}{\field{C}}

\newcommand{\R}{\field{R}}
\newcommand{\Z}{\field{Z}}
\newcommand{\PO}{\mathcal{P}}

\newcommand{\SigmaM}{\Sigma(M)}

\newcommand{\DeltaM}{\Delta(M)}

 \DeclareMathOperator{\Hom}{Hom}
\DeclareMathOperator{\rank}{rank} \DeclareMathOperator{\tr}{tr}
\DeclareMathOperator{\DHF}{DH} 
 
\DeclareMathOperator{\Sign}{Sign}

\DeclareMathOperator{\ch}{ch}
\DeclareMathOperator{\ind}{ind}

\newcommand{\lam}{\lambda}
\newcommand{\var}{\varphi}
\newcommand{\varv}{\varphi^v}

\newcommand{\hatvar}{\hat{\varphi}}
\newcommand{\hatvarv}{\hat{\varphi}^v}

\newcommand{\hatM}{\hat{M}}
\newcommand{\hatV}{\hat{V}_x}
\newcommand{\hatU}{\hat{U}_x}

\newcommand{\hatF}{\hat{F}}
\newcommand{\hatH}{\hat{H}}

\newcommand{\sigmn}{\Sigma^{(n)}}
\newcommand{\sigmk}{\Sigma^{(k)}}
\newcommand{\sigmone}{\Sigma^{(1)}}
\newcommand{\uv}{\langle u_i^I,v\rangle}

\def\l{\langle}
\def\r{\rangle}

\newcommand{\Proj}{\field{P}}

\newcommand{\img}{\sqrt{-1}}

\title{Elliptic genera, torus orbifolds and multi-fans}
\author{Akio Hattori and Mikiya Masuda}
\address{Graduate School of Mathematical Science, University of Tokyo,
Tokyo, Japan;
Department of Mathematics, Osaka city University, Osaka, Japan}
\email{hattori@ms.u-tokyo.ac.jp; masuda@sci.osaka-cu.ac.jp}
 
\begin{document}
\maketitle


\begin{abstract}
Multi-fan is an analogous notion of fan in toric theory. 
Fan is a combinatorial 
object associated to a toric variety. Multi-fan is associated to 
an orbifold with an action of half the dimension of the orbifold. 
In this paper the equivariant elliptic genus and 
the equivariant orbifold elliptic genus of multi-fans are defined and their 
character formulas are exhibited. A vanishing theorem concerning 
elliptic genus of multi-fans of global type 
and its applications to toric varieties are given. 
\end{abstract}



\section{Introduction}
\label{sec:intro}
A torus orbifold is an oriented closed orbifold of even
dimension which admits an action of a torus of half the dimension 
of the orbifold with some orientation data concerning codimension 
two fixed point set components of circle subgroups and with some
restrictions on isotropy groups of points of the orbifold. Typical examples 
are complete toric varieties with simplicial fan. To a toric variety
there corresponds a fan, and that correspondence is ono-to-one.
In particular algebro-geometric properties of a toric variety
can be described in terms of combinatorial properties of the
corresponding fan in principle, see e.g. \cite{O}. In a similar way,
to a torus orbifold there corresponds a multi-fan, a generalaization
of the notion of fan. The notion of multi-fan 
was introduced in \cite{M} and a combinatorial
theory of multi-fans was developed in \cite{HM}. In particular it
was shown there that every complete simplicial multi-fan of dimension greater
than two can be realized as the one associated with a torus orbifold.
It should be noticed that a torus orbifold also determines a set of 
vectors which generates the one dimesnional cones of the associated
multi-fan. It turns out that 
many topological invariants of a torus orbifold can be described in
terms of the multi-fan and the set of generating vectors associated 
with it. 

The purpose of the present note is to discuss elliptic genera for
torus orbifolds and multi-fans. Two sorts of elliptic genus are defined 
for stably almost complex orbifolds. One is the direct 
generalization of elliptic genus for 
stably almost complex manifolds to stably almost complex 
orbifolds which we shall denote by $\var$. 
The other, denoted by $\hatvar$, is the so-called orbifold elliptic 
genus which 
has its origin in string theory. Correspondingly we can define two sorts
of elliptic genus for a pair of complete simplicial multi-fan and generating
vectors. 
Note that, though there may not be a stably almost complex complex 
structure on a torus orbifold in general, one can still define elliptic
genera $\var$ and $\hatvar$ for pairs of complete simplicial multi-fan 
and generating vectors, and 
therefore define them for general torus orbifolds vice versa.

Borisov and Libgober gave a beautiful formula for elliptic genus 
$\var$ of complete non-singular toric varieties in \cite{BL1}.
Theorem \ref{thm:BL} and Theorem \ref{thm:BL2} describe 
similar formulae expressing the equivariant
elliptic genera $\var$ and $\hatvar$ of a complete simplicial
multi-fan as a virtual character of 
the associated torus.
The starting point of \cite{BL1} was the sheaf cohomology of toric
varieties. Our starting point is the fixed point formula of the
Atiyah-Singer type due to Vergne applied to the action of the torus. 

A remarkable feature of elliptic genera is their rigidity property.
If the circle group acts on a closed almost complex (or
more generally stably almost complex) manifold whose first Chern class is divisible
by a positive integer $N$ greater than $1$, then its equivariant
elliptic genus of level $N$ is rigid, that is, it is a constant 
character of the circle group. It was conjectured by Witten \cite{W}
and proved by Taubes \cite{T}, Bott-Taubes \cite{BT} and Hirzebruch 
\cite{H}. Liu \cite{L}
found a simple proof using the modularity of elliptic genera. 
Applying this to a non-singular complete toric variety we see
that its elliptic genera $\var$ with $(-y)^N=1,\ -y\not=1$ 
of level $N$ is rigid if its first Chern
class is divisible by $N$. Moreover, using a vanishing theorem due to
Hirzebruch \cite{H}, we can show that the genus actually
vanishes. In this note we shall extend this result to global torus orbifolds
and the corresponding type of multi-fans and generating vectors. 
We call this particular type the global type. Here a torus orbifold $M$ is 
called a global torus orbifold if $M$ is the quotient of a torus manifold 
$\tilde{M}$ by a finite subgroup of the torus acting on $\tilde{M}$.
The rigidity of elliptic genus and orbifold elliptic genus for general
orbifolds do not hold in general. We shall return to this point and
related topics in another paper.

The fact that $\var$ and $\hatvar$ are virtual characters of a torus 
is proved by using the 
multiplicity formula for Duistermaat-Heckman function for 
multi-polytopes given in \cite{HM}. The Chern class
is also defined for a pair of a multi-fan and its generating vectors. 
The rigidity and vanishing property of
elliptic genus of level $N$ can be 
formulated for a pair of multi-fan and its generating vectors of 
the global type. 
One of the main results is Corollary \ref{coro:vanishing} which 
states that, if the first Chern class of a pair of a multi-fan and 
its genrating vectors of the global type is 
divisible by $N$, then its elliptic genus 
of level $N$ vanishes. The proof 
of rigidity and vanishing property follows the 
idea of the proof given in \cite{H} translated in combinatorial terms.  
When $N=2$, 
the torus manifold is a spin manifold. The corresponding multi-fan
might be called a non-singular spin multi-fan. As a corollary we see that its
signature vanishes in this case.

The equivariant $T_y$-genus can be considered as a special value 
of equivariant elliptic genus. It was shown in \cite{HM} that the 
equivariant $T_y$-genus
of a complete multi-fan was rigid, and a formula
for the $T_y$-genus was given.
We shall give another proof of
this formula using the character formula mentioned before.  
Moreover, if the first Chern class is divisible by $N$,
then $T_y$-genus vanishes for $(-y)^N=1$. 
One can derive some
applications from this fact. For example if $\Delta$ is a complete 
non-singular multi-fan of dimension $n$ 
with first Chern class $c_1(\Delta)$ divisible by $N$ and with
non-vanishing Todd genus, then $N$ 
must be equal to or less than $n+1$ (Proposition \ref{prop:N=n+1}). 
In the extremal case $N=n+1$, if $\Delta$ is assumed to be a 
complete non-singular ordinary fan, then $\Delta$ must be isomorphic
to the fan of projective space $\Proj^n$. 
Hence a complete non-singular toric variety $M$ of dimension $n$ with 
$c_1(M)$ divisible by $n+1$ must be isomorphic to
$\Proj^n$ as toric variety (Corollary \ref{coro:projective space}). 
We show furthermore that,  
in case $c_1(M)$ is divisible by $n$, $M$ is isomorphic to a certain 
projective space bundle over $\Proj^1$ (Corollary \ref{coro:projective 
bundle}). The authors are grateful to T. Oda and T. Fujita for informing them 
that these results can be obtained by standard arguments 
in algebraic geometry at least for projective toric varieties. They 
are also grateful to O. Fujino who communicated to them his proof of 
these results including the case of singular varieties \cite{Fuj}. 

The paper is organized as follows. In Section 2 we recall some basic facts
about multi-fans from \cite{HM}. 
In Section 3 we define the elliptic genus and orbifold elliptic 
genus of a pair of a multi-fan and its generating vectors and derive 
the character formulae (Theorem \ref{thm:BL} and \ref{thm:BL2}). 
A formula for the $T_y$-genus is also given.
In Section 4 we define equivariant
first Chern class of a pair of a multi-fan and its generating vectors 
and discuss properties concerning
its divisibility. In Section 5 the proof of the rigidity
of elliptic genus of level $N$ for the global type is given. The
main results are Theorem \ref{theo:rigid} and Corollary 
\ref{coro:vanishing}. 
Section 6 is devoted to applications. 
In the last section we shall recall the fixed point formula, 
the formulae for
elliptic genus and orbifold elliptic genus for almost complex orbifolds, 
and give the explicit
formulae for elliptic genus and orbifold elliptic genus of torus
orbifolds from which the corresponding formulae for multi-fans 
are deduced.


\section{Multi-fans}\label{sec:multi}
We refer to \cite{HM} for notions and notations
concerning multi-fans and torus orbifolds. We shall summarize some of 
them in the sequel. 
Let $L$ be a lattice of rank n and $\Delta=(\Sigma,C,w^{\pm})$ 
an $n$-dimensional 
\emph{simplicial multi-fan} in $L$ (the notation $N$ was used in \cite{HM}
instead of $L$). Here $\Sigma$ is an
augmented simplicial set, that is, 
$\Sigma$ is a simplicial set with empty set 
$*=\emptyset$ added as $(-1)$-dimenional simplex.
$\sigmk$ denotes the $k-1$ skeleton of $\Sigma$ so that 
$*\in \Sigma^{(0)}$. 
We assume that $\Sigma=\sum_{k=0}^n\sigmk$, and $\sigmn \not=\emptyset$.
$C$ is a map from $\sigmk$ into the set of $k$-dimensional 
strongly convex rational
polyhedral cones in the vector space $L_\R=L\otimes \R$ for 
each $k$ such that, if $J$ is a face of $I$, then $C(J)$ is
a face of $C(I)$.  

$w^{\pm}$ are maps $\sigmn \to \Z_{\ge 0}$ which, when $\Sigma$
is complete, satisfy certain compatibility conditions, as
we shall explain below. 
We set $w(I)=w^+(I)-w^-(I)$.
A vector $v\in L_\R$ will be called \emph{generic} if $v$ does
not lie on any linear subspace spanned by a cone in 
$C(\Sigma)$ of dimesnsion less than $n$. For a generic
vector $v$ we set $d_v=\sum_{v\in C(I)}w(I)$, where
the sum is understood to be zero if there is no such $I$.
We call a multi-fan $\Delta=(\Sigma,C,w^\pm)$ of dimension $n$ 
{\it pre-complete} if 
the integer $d_v$ is independent of
the choice of generic vectors $v$. We call this integer
the {\it degree} of $\Delta$ and denote it by $\deg(\Delta)$.

For each $K\in \Sigma$ we set
\[ \Sigma_K=\{ J\in\Sigma\mid K\subset J\}.\]
It inherits the partial ordering from $\Sigma$ and becomes an augmented
simplicial set where $K$ is the unique 
minimum element in $\Sigma_K$. 
Let $(L_K)_\R$ be the linear subspace of $L_\R$ generated by $C(K)$.
Let $L^K_\R$ be the quotient space of $L_\R$ by $(L_K)_\R$ and 
$L^K$ the image of $L$ in $L^K_\R$. $L^K_\R$ is identified with
$L^K\otimes \R$. 
For $J\in \Sigma_K$ we define $C_K(J)$ to be
the cone $C(J)$ projected on $L^K_\R$. 
We define two functions 
\[ {w_K}^\pm\colon  \Sigma_K^{(n-|K|)}\subset \Sigma^{(n)} \to \Z_{\ge 0}\]
to be the restrictions of $w^\pm$ to $\Sigma_K^{(n-|K|)}$. The triple 
$\Delta_K:=(\Sigma_K,C_K,{w_K}^\pm)$ is a multi-fan in $L^K$ 
and is called the  
\emph{projected multi-fan} with respect to $K\in \Sigma$. 
If $K=\emptyset$ then $\Delta_K =\Delta$. 
A pre-complete multi-fan $\Delta=(\Sigma,C,w^\pm)$ is 
said to be 
\emph{complete} if the projected multi-fan $\Delta_K$ is pre-complete for any 
$K\in \Sigma$. A multi-fan is complete if and only if the projected
multi-fan $\Delta_J$ is pre-complete for any $J\in \Sigma^{(n-1)}$.

Let $M$ be an oriented closed manifold of dimension $2n$ with an 
effective action of an $n$-dimensional torus $T$. We assume further 
that the fixed point set $M^T$ is not empty. There is a 
finite number of subcircles of $T$ such that the fixed point
set of each subcircle has codimension 2 components. Let $\{M_i\}_{i=1}^r$
be those components which have non-empty intersection with $M^T$.
We call $M$ a \emph{torus manifold} if a preferred orientation of 
each $M_i$ is given. The $M_i$ are called \emph{characteristic submanifolds}. 
A multi-fan 
$\Delta(M)=(\Sigma(M),C(M),w^\pm(M))$
in the lattice $H_2(BT)$ is associated with $M$, where
$BT$ is the classifying space of $T$ 
(homology is taken  
with coefficients in the integers, unless otherwise specified). 
The (augmented) simplicial set $\Sigma(M)$ is defined by 
\[ \Sigma(M)=\{I\subset \{1,\ldots ,r\}| M_I=(\cap_{i\in I}M_i)^T\not=\emptyset \}. \]
We make convention that $M_I=M$ for $I=*=\emptyset$. 
The cones $C(M)(I)$ are defined as follows.
Let $\nu_i$ denote the normal bundle of $M_i$ in $M$. It is an oriented
$2$-plane bundle with the orientation induced by those of $M_i$ 
and $M$, and, as such, it is regarded as a complex line bundle. If $S_i$
is the subcircle which fixes $M_i$ pointwise, then $S_i$ acts effectively
on each fiber of $\nu_i$ as complex automorphism. Hence there is a unique 
isomorphism
$\rho :S^1 \to S_i$ such that $\rho_i(t)$ acts as complex multiplication
by $t$. Thus $\rho_i$ defines a primitive element $v_i\in \Hom(S^1,T)$.
We identify 
$\Hom(S^1,T)$ with $L=H_2(BT)$. If $I$ is in $\sigmn$, 
then $\{v_i\}_{i\in I}$ is a basis of $L$.
If $I$ is in $\Sigma$, the cone $C(M)(I)$ is defined to be the cone 
generated by $\{v_i\}_{i\in I}$ in $L\otimes \R$.  
The fixed point set $M^T$ coincides with the union $\cup_{I\in \Sigma(M)^{(n)}}M_I$.
For each $p\in M^T$ let $\epsilon_p=\pm 1$ be the ratio of two orientations
at $p$, one induced from that of $M$ and the other determined as the 
intersection of oriented submanifolds $\{M_i\}_{i\in I}$. The number
$w(M)^+(I)$ (respectively $w(M)^-(I)$) is defined as the number of 
$p\in M_I$ with $\epsilon_p=+1$ (respectively $\epsilon_p=-1$).   
$\Delta(M)$ is a complete multi-fan. 
If $K\in \Sigma(M)$ then the projected multi-fan $\Delta(M)_K$ is closely 
related to the multi-fan associated with $M_K=\cap_{i\in K}M_i$, where
$M_K$ is regarded as a union of torus manifolds, see \cite{HM} for details.

Let $\Delta=(\Sigma,C,w^\pm)$ be a multi-fan in $L$. 
If $T$ denotes the torus $L_\R/L$, then $L$ can be canonically
identified with $H_2(BT)$.
Then there is a unique
primitive vector $v_i\in L=H_2(BT)$ which generates the cone $C(i)$
for each $i\in \sigmone$. $\Delta$ is called non-singular if
$\{v_i\mid i\in I\}$ is a basis of the lattice $L=H_2(BT)$ for each 
$I\in \sigmn$. Thus the multi-fan associated with a torus manifold
is a complete non-singular multi-fan.

It is sometimes more convenient to consider a set of vectors
${\mathcal V}=\{v_i\in L\}_{i\in \sigmone}$ such that each $v_i$ generates
the cone $C(i)$ in $L_\R$ but is not necessarily primitive. 
This is the case for multi-fans associated with 
torus orbifolds. A torus orbifold is a closed
oriented orbifold of even dimension with an effective action of 
a torus of half the dimension of the orbifold with some additional
condition. We refer to \cite{HM} for details. 
A set of codimension 2 suborbifolds $M_i$ called \emph{characteristic 
suborbifolds} is similarly defined
as in the case of torus manifolds. To each subcircle $S_i$
which fixes $M_i$ pointwise there is some finite cover 
$\tilde{S}_i$ and an effective action of $\tilde{S}_i$ on the 
orbifold cover of each fiber of the normal bundle.
This defines a vector $v_i$ in $\Hom(S^1,T)=H_2(BT)=L$ as before.
In this way a multi-fan $\Delta(M)$ and a set of vectors
${\mathcal V}(M)=\{v_i\}_{i\in \sigmone}$ are associated to
a torus orbifold $M$. 

Hereafter multi-fans are assumed to be complete and simplicial, 
and a set of vectors
${\mathcal V}=\{v_i\in L\}_{i\in \sigmone}$ as above is 
associated to each multi-fan $\Delta=(\Sigma,C,w^\pm)$. 
In case $\Delta$ is non-singular it is further assumed
that all the $v_i$ are primitive.
If $I$ is in $\sigmn$, then $\{v_i\}_{i\in I}$ becomes
a basis of vector space $L_\R$.
In case $\Delta$ is non-singular it is a basis of the lattice
$L$. In general, for $I\in \sigmn$, we define $L_{I,\V}$ to be 
the sublattice of $L$
generated by $\{v_i\}_{i\in I}$. 

Let $L_{I,\V}^*$ be the dual lattice of
$L_{I,\V}$ and and $\{u_i^I\}$ the basis of $L_{I,\V}^*$ dual 
to $\{v_i\}_{i\in I}$. 
We identify $L_{I,\V}^*$ with the lattice in $L_\R^*$
given by 
\[ \{u\in L_\R^*\mid \langle u,v\rangle \in \Z,
   \ \ \text{for any}\ \  v\in L_{I,\V}\},\] 
where $\l u,v\r$ is the dual pairing. 
For $h\in L/L_{I,\V}$ and $u\in L_{I,\V}^*$ we define 
\[ \chi_I(u,h)=e^{2\pi\sqrt{-1}\langle u,v(h)\rangle},\]
where $v(h)\in L$ is a representative of $h$. If one fixes $u$,
$h \mapsto \chi_I(u,h)$ gives a character of the group $L/L_{I,\V}$. 

The dual lattice
$L^*=H^2(BT)\subset H^2(BT;\R)$ is canonically identified 
with $\Hom(T,S^1)$. The latter is embedded in the character ring
$R(T)$. In fact $R(T)$ can be considered as the group ring $\Z[L^*]$ of 
the group $L^*=\Hom(T,S^1)$. It is convenient to write the element 
in $R(T)$ corresponding to $u\in H^2(BT)$ by $t^u$. The 
homomorphism $v^*:R(T)\to R(S^1)=\Z[t,t^{-1}]$ induced by an element 
$v\in H_2(BT)=\Hom(S^1,T)$ can be written in the form
\[ v^*(t^u)=t^{\langle u,v\rangle}, \]
where $t^m\in R(S^1)$ is such that $t^m(g)=g^m$ for $g\in S^1$.

More generally, set $L_\V=\bigcap_{I\in \sigmn}L_{I,\V}$, and 
let $L_{\V}^*$ be the dual lattice of $L_\V$. 
$L_{\V}^*$ contains all $L_{I,\V}^*$ and is generated by all the 
$u_i^I$'s. The group ring $\Z[L_{\V}^*]$ 
contains $\Z[L^*]=R(T)$ and has a basis $\{t^u|u\in L_{\V}^*\}$ with 
multiplication determined by the addition in $L_{\V}^*$:
\[  t^ut^{u'}=t^{u+u'}. \]
If $v$ is a vector in $L_\V$, then $v$ determines a homomorphism
$v^*:\Z[L_{\V}^*]\to R(S^1)=\Z[t,t^{-1}]$ sending $t^u$ to
$t^{\langle u,v\rangle}$. If we vary $v$ then 
$v^*(t^u)$ determines $t^u$. 

 Similarly if $v_1$ and $v_2$ are vectors 
in $L$, then they define a homomorphism from a $2$-dimensional 
torus $T^2$ into $T$ and induce a homomorphism 
$(v_1,v_2)^*:\Z[L^*]\to R(T^2)=\Z[t_1,t_1^{-1},t_2,t_2^{-1}]$ 
defined by 
\[ (v_1,v_2)^*(t^u)=t_1^{\l u,v_1\r}t_2^{\l u,v_2\r}. \]
Moreover if $v_1$ and $v_2$ belong to $L_\V$, then $(v_1,v_2)^*$ 
extends to a homomorphism $\Z[L_\V^*]\to R(T^2)$.

We define the equivariant cohomology $H_T^*(\Delta)$ of a complete 
multi-fan $\Delta$ as the face ring of the simplicial complex $\Sigma$.
Namely let $\{x_i\}$ be indeterminates indexed by $\sigmone$, and 
let $R$ be the polynomial ring over the integers generated 
by $\{x_i\}$. We denote by $\mathcal{I}$ the ideal in $R$
generated by monomials
$\prod_{i\in J} x_i$ such 
that $J\notin \Sigma$. $H_T^*(\Delta)$ is by definition
the quotient $R/\mathcal{I}$. 
We regard $H^2(BT)$ as a submodule of $H_T^2(\Delta)$ by the formula
\begin{equation}\label{eq:structure}
 u=\sum_{i\in \sigmone}\langle u,v_i\rangle x_i .
\end{equation} 
This determines an $H^*(BT)$-module structure of $H_T^*(\Delta)$.
It should be noticed that this module structure depends on
the choice of vectors ${\mathcal V}$ as above.
Let $S$ be the subset of $H^*(BT)$ multiplicatively generated by
non-zero elements in $H^2(BT)$. 
If $M$ is a torus manifold, then $H_T^*(\Delta(M))$ can be embedded in
$H_T^*(M)$ divided by $S$-torsions and coincides 
with it provided some additional conditions are satisfied.

For each $I\in\sigmn$ we define the restriction homomorphism
$\iota_I^* :H_T^2(\Delta)\to L_{\V}^*$ by
\begin{equation*}\label{eq:restrict}
 \iota_I^*(x_i)=\begin{cases}
            u_i^I & \text{for}\ i\in I\\
            0     & \text{for}\ i\notin I.
            \end{cases}
\end{equation*}            
It follows from (\ref{eq:structure}) that $\iota_I^*| H^2(BT)$ is 
the identity map for any $I$, and $\sum_{I\in \sigmn}\iota_I^*$ is 
injective. Note that, if $\Delta$ is non-singular, then $\iota_I^*$ maps
$H_T^2(\Delta)$ into $H^2(BT)$.

\begin{lemm}\label{lemm:multiplicity bis}
For any $x=\sum_{i\in \sigmone}c_ix_i\in H_T^2(\Delta),\ c_i\in \Z$,
the element 
\[ \sum_{I\in\sigmn}\frac{w(I)}{|L/L_{I,\V}|}\sum_{h\in L/L_{I,\V}}
           \frac{\chi_I(\iota_I^*(x),h)t^{\iota_I^*(x)}}
         {\prod_{i\in I}(1-\chi_I(u_i^I,h)^{-1}t^{-u_i^I})}\]
in $\C[L_{\V}^*]$ actually belongs to $R(T)$.
\end{lemm}
This was proved in Corollary 7.4 of \cite{HM} with a further assumption
that $\iota_I^*(x)\in H^2(BT)$. The general case can be proved in a
similar way. The formula was also given in Corollary 12.10 of 
\cite{HM} when $\Delta$ is the multi-fan associated with a torus 
orbifold. 

We also use an extended version of Corollary 7.4 in \cite{HM}. 
Let $K\in \sigmk$ and let $\Delta_K=(\Sigma_K,C_K,w_K^{\pm})$ be 
the projected multi-fan.
If $I\in\Sigma^{(l)}$ contains $K$, then $I$ is considered 
as lying in $\Sigma_K^{(l-k)}$. In order to avoid some notational confusions
we introduce the link $\Sigma'_K$ of $K$ in $\Sigma$.
It is a simplicial set consisiting simplices $J$ such that 
$K\cup J\in \Sigma$ and $K\cap J=\emptyset$. There is an isomorphism 
from $\Sigma'_K$ to $\Sigma_K$ sending $J\in {\Sigma'_K}^{(l)}$ to 
$K\cup J\in \Sigma_K^{(l)}$. Let $K\ast\Sigma'_K$ be the join of 
$K$ (regarded as a simplicial set) and $\Sigma'_K$. Its simplices are
of the form $J_1\cup J_2$ with $J_1\subset K$ and $J_2\in \Sigma'_K$. 
The torus $T^K$ corresponding to $\Delta_K$ 
is a quotient of $T$. We consider the polynomial 
ring $R_K$ generated by $\{x_i\mid i\in K\cup{\Sigma'_K}^{(1)}\}$ 
and the ideal $\mathcal{I}_K$ generated by monomials
$\prod_{i\in J} x_i$ such 
that $J\notin K\ast\Sigma'_K$. We define the equivariant cohomology 
$H_T^*(\Delta_K)$ of
$\Delta_K$ with respect to the torus $T$ as the quotient ring
$R_K/\mathcal{I}_K$. Note that $H_T^*(\Delta_K)$ is defferent from 
$H_{T^K}^*(\Delta_K)$. 

$H^2(BT)$ is regarded as a submodule of
$H_T^2(\Delta_K)$ by a formula similar to (\ref{eq:structure}).
This defines
an $H^2(BT)$-module structure on $H_T^*(\Delta_K)$ .
The projection $H_T^2(\Delta)\to H_T^2(\Delta_K)$ 
is defined by sending $x_i$ to $x_i$ for $i\in K\cup{\Sigma'_K}^{(1)}$ and
putting $x_i=0$ for $i\notin K\cup{\Sigma'_K}^{(1)}$.  
The restriction homomorphism
$\iota_I^*:H_T^2(\Delta_K) \to L_{\V}^*$ is also defined for 
$I\in \Sigma_K^{(n-k)}$ by $\iota_I^*(x_i)=u_i^I$. If $M$ is a torus orbifold, then 
$H_T^*(\Sigma(M)_K)$ is related to the equivariant cohomology
$H_T^*(M_K)$ with respect to the group $T$ (not with respect to $T_K$),
cf. Remark \ref{rem:face ring} in Section 4. 

Given   
$x=\sum_{i\in K\cup{\Sigma'_K}^{(1)}}c_ix_i\in H_T^2(\Delta_K)\otimes\R,
\ c_i\in \R$, 
let $A^*$ be 
the affine subspace in the dual space $L_\R^*$ defined by 
$\langle u,v_i\rangle =c_i$ for $i\in K$. Then we introduce 
a collection $\mathcal{F}_K=\{F_i\mid i\in{\Sigma'}_K^{(1)}\}$
of affine hyperplanes in $A^*$ by setting 
\[ F_i= \{u\mid u\in A^*,\ \langle u,v_i \rangle =c_i \}. \]
The pair $\PO_K=(\Delta_K,\mathcal{F}_K)$ will be called a 
multi-polytope associated with $x$; see \cite{HM} for the case 
$K=\emptyset$. For $I\in \Sigma_K^{(n-k)}$, i.e. 
$I \in \sigmn$ with $I\supset K$, we put 
$u_I=\cap_{i\in I\setminus K}F_i\in A^*$. 
Note that $u_I$ is equal to $\iota_I^*(x)$. The dual vector space 
$(L^K_\R)^*$ of $L^K_\R$ is canonically identified with 
the subspace $\{u\mid \langle u,v_i\rangle =0,\ i\in K\}$ of 
$L^*_\R=H^2(BT;\R)$. It is parallel to $A^*$, and $u_i^I$ lies 
in $(L^K_\R)^*$ for $I\in \Sigma_K^{(n-k)}$ and $i\in I\setminus K$.
A vector $v\in L^K_\R$ is called generic if
$\langle u_i^I,v\rangle \not=0$ for any $I\in \Sigma_K^{(n-k)}$ and
$i\in I\setminus K$. 
The image in $L_\R^K$ of a generic vector in $L_\R$ is generic. 
We take a generic vector $v\in L^K_\R$, and define, for 
$I\in \Sigma_K^{(n-k)}$ and $i\in I\setminus K$,
\[ (-1)^I:=(-1)^{\#\{j\in I\setminus K\mid \langle u_j^I,v\rangle >0\}}
\quad \text{and}\quad
(u_i^I)^+ :=
      \begin{cases}
      u_i^I & \text{if}\  \langle u_i^I,v\rangle >0\\
     -u_i^I & \text{if}\  \langle u_i^I,v\rangle <0. 
      \end{cases}            \]
We denote by $C_K^*(I)^+$ the cone in $A^*$ spanned by the
$(u_i^I)^+,\ i\in I\setminus K,$ with apex at $u_I$, and by $\phi_I$ its
characteristic function. With these understood, we define a
function $\DHF_{\PO_K}$ on $A^*\setminus \cup_iF_i$ by
\[ \DHF_{\PO_K} :=\sum_{I\in \Sigma_K^{(n-k)}}(-1)^Iw(I)\phi_I. \]
As in \cite{HM} we call this function the \emph{Duistermaat-Heckman 
function} associated with $\PO_K$. 

\begin{lemm}
The support of the function $\DHF_{\PO_K}$ is bounded, and
the function is independent of the choice of generic vector $v$.
\end{lemm}
The proof is similar to that of Lemma 5.4 in \cite{HM}.
We shall denote by $\PO_{K+}$
the multi-polytope associated with 
$x_+=\sum_{i\in K}c_ix_i +\sum_{i\in{\Sigma'_K}^{(1)}}(c_i+\epsilon)x_i$ 
where $0<\epsilon<1$. The following theorem is a generalization 
of Corollary 7.4 in \cite{HM}. 

\begin{theo}\label{theo:multiplicity}
Let $\Delta$ be a complete simplicial multi-fan.
Let $x=\sum_{i\in K\cup{\Sigma'_K}^{(1)}}c_ix_i\in H_T^2(\Delta_K)$ be as 
above with all $c_i$ integers,
and let $\PO_{K+}$ be defined as above. Then
\begin{equation*}
\sum_{u\in A^*\cap L^*}\DHF_{\PO_{K+}}(u)t^u
 =\sum_{I\in \Sigma_K^{(n-k)}}\frac{w(I)}{|L/L_{I,\V}|}\sum_{h\in L/L_{I,\V}}
 \frac{\chi_I(\iota_I^*(x),h)t^{\iota_I^*(x)}}
 {\prod_{i\in I\setminus K}(1-\chi_I(u_i^I,h)^{-1}t^{-u_i^I})}.
\end{equation*}
In particular the right hand side belongs to $R(T)$.
\end{theo}
The proof is similar to that of Corollary 7.4 in \cite{HM}.
Applying $(v_1,v_2)^*$ to the both sides of the above equality  
we get
\begin{coro}\label{coro:multiplicity}
Let $v_1$ and $v_2$ be generic vectors in $L$ such that 
$\l\u,v_1\r$ and $\l\u,v_2\r$ are integers for all $I\in\Sigma_K^{(n-k)}$ 
and let $x=\sum_{i\in K\cup{\Sigma'_K}^{(1)}}c_ix_i\in H_T^2(\Delta_K)$ 
with all $c_i$ integers. Then 
\begin{equation*}
\sum_{I\in \Sigma_K^{(n-k)}}\frac{w(I)}{|L/L_{I,\V}|}\sum_{h\in L/L_{I,\V}}
 \frac{\chi_I(\iota_I^*(x),h)
  t_1^{\l\iota_I^*(x),v_1\r}t_2^{\l\iota_I^*(x),v_2\r}}
   {\prod_{i\in I\setminus K}
  (1-\chi_I(u_i^I,h)^{-1}t_1^{\l-u_i^I,v_1\r}t_2^{\l-u_i^I,v_2\r})} 
\end{equation*}
belongs to $R(T^2)=\Z[t_1,t_1^{-1},t_2,t_2^{-1}]$.
\end{coro}

For $I\in \Sigma_K^{(n-k)}$ let $G_I$ be the 
subgroup of the permutation group of 
$I$ consisting of those elements which are identity on $K$. Let 
$\mathcal{L}_I$ be the set of all linear forms 
$\sum_{i\in I}m_iu_i^I$ with integer coefficients $m_i$. 
The group $G_I$ acts on $\mathcal{L}_I$. Let $\mathcal{O}_I$
denote the set of orbits of that action. If $I'$ is also in
$\Sigma_K^{(n-k)}$, take a bijection $f:I\to I'$ which is the
identity on $I\cap I'$. It induces a bijection 
$f_*:\mathcal{L}_I\to \mathcal{L}_{I'}$ defined by 
$f_*(u_i^I)=u_{f(i)}^{I'}$. This in turn induces a bijection 
$f_*:\mathcal{O}_I\to \mathcal{O}_{I'}$ . It is easy to see that 
the latter $f_*$ does not depend on the choice of particular $f$. 
Thus we shall write it $\tau_{I,I'}$. 

The following lemma will be useful later. 
\begin{lemm}\label{lemm:integrality}
Let $v_1$ and $v_2$ be generic vectors in $L$ as in Corollary 
\ref{coro:multiplicity} and 
let $\mathcal{A}=\{\alpha_I\in \mathcal{O}_I\}_{I\in \Sigma_K^{(n-k)}}$ be 
a collection which satisfies the relations 
\[ \tau_{I,I'}(\alpha_I)=\alpha_{I'}\ \ 
\text{for any $I,I'\in\Sigma_K^{(n-k)}$}. 
\]
Then the expression
\[ \sum_{I\in\Sigma_K^{(n-k)}}\frac{w(I)}{|L/L_{I,\V}|}
 \sum_{l\in \alpha_I}\sum_{h\in L/L_{I,\V}}
 \frac{\chi_I(l,h)t_1^{\l l,v_1\r}t_2^{\l l,v_2\r}}
 {\prod_{i\in I\setminus K}(1-\chi_I(u_i^I,h)^{-1}
 t_1^{\l -u_i^I,v_1\r}t_2^{\l -u_i^I,v_2\r})}\]
belongs to $\Z[t_1,t_1^{-1},t_2,t_2^{-1}]$.
\end{lemm}
\begin{proof}
We define the support $supp(l)$ of a linear form 
$l=\sum_{i\in I}m_iu_i^I \in \mathcal{L}_I$
to be the set $\{i\in I\vert m_i\not=0\}$. The cardinal number 
of $supp(l)$ is called the length of $l$ and is denoted by 
$|l|$. 
The length is invariant under the action of $G_I$,
so that the length $|\alpha_I|$ of $\alpha_I\in \mathcal{O}_I$
is defined as that of a linear form contained in $\alpha_I$. 
The bijections $\tau_{I,I'}$ preserve length.
The proof will proceed by induction on the common length of 
the $\alpha_I$ in the hypothesis of Lemma. 

If $|\alpha_I|=0$, then $\alpha_I$ consists of $0$.
Applying Corollary 
\ref{coro:multiplicity} to $x=0$ we see that the statement of 
Corollary is true in this case. 

Suppose that the 
statement is true for all $\alpha_I$ of length less than $r>0$.
If $\mathcal{L}$ denotes the set of all linear forms 
$\sum_{i\in K\cup{\Sigma'_K}^{(1)}}m_ix_i$, there is 
a natural injection $j^I_*: \mathcal{L}_I\to\mathcal{L}$ 
sending $u_i^I$ to $x_i$. Then it is clear that
$ \iota_I^*(j^I_*(l))=l $ 
for any $l\in \mathcal{L}_I$. On the other hand, 
$\iota_I^*(j^{I'}_*(f_*(l)))$ has length less than that of $l$ 
for any bijection $f:I\to I'$ fixing elements of $I\cap I'$ 
unless $supp(l)\subset I\cap I'$. In the latter case we have
\[ j^{I'}_*(f_*(l))=j^I_*(l). \]

Let $\mathcal{L}(\mathcal{A})$ be the totality of linear forms
of the form $j^I_*(l)$ with $l\in \alpha_I$ and 
$\alpha_I\in \mathcal{A}$. Then, for each $l\in \mathcal{L}(\mathcal{A})$ 
and $I$, 
the linear form $\iota_I^*(l)$ either belongs to $\alpha_I$ 
or has length less than $r$. Moreover each $l_1\in \alpha_I$ 
appears exactly once in this way. It is also easy to see that 
the totality of transforms of $\iota_I^*(l)$ with length 
less than $r$ by permutations of $I$ fills some orbits in 
$\mathcal{O}_I$ with multiplicities. If this set (with 
multiplicity) is denoted by $\mathcal{L}_I(\mathcal{A})$, 
then one has
\[ \begin{split} 
  &\sum_{I\in\Sigma_K^{(n-k)}}\sum_{l\in \alpha_I}
 \frac{w(I)}{|L/L_{I,\V}|}\sum_{h\in L/L_{I,\V}}
 \frac{\chi_I(l,h)t_1^{\l l,v_1\r}t_2^{\l l,v_2\r}}
 {\prod_{i\in I\setminus K}(1-\chi_I(u_i^I,h)^{-1}
  t_1^{\l -u_i^I,v_1\r}t_2^{\l -u_i^I,v_2\r})} \\
  =& \sum_{l\in\mathcal{L}(\mathcal{A})}\sum_{I\in\Sigma_K^{(n-k)}}
  \frac{w(I)}{|L/L_{I,\V}|}\sum_{h\in L/L_{I,\V}}
  \frac{\chi_I(\iota_I^*(l),h)t_1^{\l \iota_I^*(l),v_1\r}
  t_2^{\l\iota_I^*(l),v_2\r}}
 {\prod_{i\in I\setminus K}(1-\chi_I(u_i^I,h)^{-1}
 t_1^{\l -u_i^I,v_1\r}t_2^{\l -u_i^I,v_2\r})}  \\
 -&\sum_{I\in\Sigma_K^{(n-k)}}\sum_{l\in \mathcal{L}_I(\mathcal{A})}
 \frac{w(I)}{|L/L_{I,\V}|}\sum_{h\in L/L_{I,\V}}
 \frac{\chi_I(l,h)t_1^{\l l,v_1\r}t_2^{\l l,v_2\r}}
 {\prod_{i\in I\setminus K}(1-\chi_I(u_i^I,h)^{-1}
 t_1^{\l -u_i^I,v_1\r}t_2^{\l -u_i^I,v_2\r})} 
 \end{split}
\]
Since the first and second terms in the right hand side belong 
to $\Z[t_1,t_1^{-1},t_2,t_2^{-1}]$ by Corollary \ref{coro:multiplicity} 
and induction assumption, the left hand side does so. 
\end{proof}


\section{Elliptic genera of multi-fans}
\label{sec:elliptic of multi}
Let $\Delta$ be a complete simplicial multi-fan in a lattice $L$ 
and $\V=\{v_i\}_{i\in\sigmone}$
a set of prescribed vectors as in Section 2. We shall define 
the (equivariant) elliptic genus $\var(\Delta,\V)$ and 
the (equivariant) orbifold elliptic genus $\hatvar(\Delta,\V)$ of 
the pair $(\Delta,\V)$. The definitions of 
$\var(\Delta,\V)$ and $\hatvar(\Delta,\V)$ are 
such that, if $M$ is an almost complex (or more
generally stably almost complex) torus orbifold, then $\var(\Delta(M),\V)$ 
and $\hatvar(\Delta(M),\V)$ coincide with those of $M$ 
expressed by the fixed point formula. 
These facts will be explained in Section 7.

We first consider the function $\Phi(z,\tau)$ of 
$z$ in $\C$ and $\tau$ in the upper half plane $\mathcal{H}$
given by the following formula.
\[ \Phi(z,\tau)=(t^{\frac12}-t^{-\frac12})\prod_{k=1}^{\infty}
 \frac{(1-tq^k)(1-t^{-1}q^k)}
      {(1-q^k)^2},  \]
where $t=e^{2\pi\img z}$ and $q=e^{2\pi\img\tau}$. Note that $|q|<1$. 
Let $A=\begin{pmatrix}a & b\\
                      c & d \end{pmatrix} \in SL_2(\Z)$, and  
put $A(z,\tau)=(\frac{z}{c\tau +d},\frac{a\tau +b}{c\tau +d})$.
$\Phi$ is a Jacobi form and satisfies the following transfomation
formulae, cf. \cite{HBJ}. 
\begin{align}
 \Phi(A(z,\tau))&=(c\tau +d)^{-1}
  e^{\frac{\pi\img cz^2}{c\tau +d}}\Phi(z,\tau), 
  \label{eq:Transformula1}\\
 \Phi(z+m\tau +n,\tau)&=(-1)^{m+n}e^{-\pi\img(2mz+m^2\tau)}\Phi(z,\tau)
  \label{eq:Transformula2}
\end{align}
where $m,n\in \Z$.

For $\sigma\in \C$ we set 
\[ \phi(z,\tau,\sigma)=\frac{\Phi(z+\sigma,\tau)}{\Phi(z,\tau)}
=\zeta^{-\frac12}\frac{1-\zeta t}{1-t}
\prod_{k=1}^{\infty}\frac{(1-\zeta tq^n)(1-\zeta^{-1}t^{-1}q^n)}
 {(1-tq^n)(1-t^{-1}q^n)}, \]
where $\zeta =e^{2\pi\img\sigma}$. 
From \eqref{eq:Transformula1} and \eqref{eq:Transformula2} 
the following transformation formulae 
for $\phi$ follow:
\begin{align}
 \phi(A(z,\tau),\sigma)&=
  e^{\pi\img c(2z\sigma +(c\tau +d)\sigma^2)}\phi(z,\tau,(c\tau +d)\sigma), 
  \label{eq:transformula1}\\
  \phi(z+m\tau +n,\tau,\sigma)&=e^{-2\pi\img m\sigma}\phi(z,\tau,\sigma)=
  \zeta^{-m}\phi(z,\tau,\sigma).
  \label{eq:transformula2}
\end{align}

In the sequel we fix the set $\V$ and put $H_I=L/L_{I,\V}$. 
Let $v\in L_{\V}$ be a generic vector. 
The (equivariant) \emph{elliptic genus} $\varv(\Delta,\V)$ 
along $v$ and 
the (equivariant) \emph{orbifold elliptic genus} 
$\hatvarv(\Delta,\V)$ along $v$
of the pair $(\Delta,\V)$ are defined by
\begin{equation}\label{eq:varv}
 \varv(\Delta,\V)=\sum_{I\in\sigmn}\frac{w(I)}{|H_I|}
  \sum_{h\in H_I}
   \prod_{i\in I}\phi(\langle u_i^I,-zv-v(h)\rangle,\tau,\sigma),
\end{equation} 
and
\begin{equation}\label{eq:hatvarv}
 \begin{split}
 \hatvarv(\Delta,\V) 
 = &\sum_{I\in\sigmn}
   \frac{w(I)}{|H_I|}\cdot \\
   &\sum_{(h_1,h_2)\in H_I\times H_I}
  \prod_{i\in I}\zeta^{\l u_i^I, v(h_1)\r}
 \phi(\l u_i^I, -zv+\tau v(h_1)-v(h_2)\r, \tau, \sigma).
\end{split} 
\end{equation}
where $v(h),v(h_1),v(h_2)\in L$ are representatives of $h,h_1,h_2\in H_I$
respectively.  
The above expressions give well-defined functions independent of 
the choice of representatives
$v(h),v(h_1),v(h_2)$ as is easily seen from \eqref{eq:transformula2}.
They are meromorphic functions in the variables $z,\tau,\sigma$ 
and sometimes written as $\varv(\Delta,\V;z,\tau,\sigma)$
and $\hatvarv(\Delta,\V;z,\tau,\sigma)$ to 
emphasize the variables. 

For each $K\in \sigmk$ with $k>0$ let $L_K$ be the kernel 
of the projection map
$L\to L^K$ and let $L_{K,\V}$ be the sublattice of $L_K$
generated by $v_i \in K$. We set 
 $H_K=L_K/L_{K,\V}$. If $J\subset K$ then we have
$ L_J\cap L_{K,\V}=L_{J,\V} $, and hence 
$H_J$ is canonically 
embedded in $H_K$. We set
\[\hat{H}_K=H_K\setminus \bigcup_{J\subsetneqq K}H_J .\]
The subset $\hat{H}_K$ is characterized by
\begin{equation}\label{eqn:0}
 \hat{H}_K=\{h\in H_K\mid \langle u_i^K,v(h)\rangle \not\in \Z 
\quad \text{for any $i\in K$} \}, 
\end{equation}
where $\{u_i^K\}$ is the basis of $L_{K,\V}^*$
dual to the basis $\{v_i\}_{i\in K}$ of $L_{K,\V}$ and $v(h)\in L_K$ is a 
representative of $h\in H_K$. 
For the minimum element $*=\emptyset\in \Sigma^{(0)}$
we set $\hat{H}_*=H_*=0$. 

If $K$ is contained in $I\in \sigmn$, then the canonical map
$L_{I,\V}^*\to L_{K,\V}^*$ sends $u_i^I$ to $u_i^K$ for $i\in K$
and to $0$ for $i\in I\setminus K$. Therefore, if $h$ is in $H_K$, 
then $\l u_i^I,v(h)\r=0$ for $i\in I\setminus K$, and 
$\l u_i^I,v(h)\r=\l u_i^K,v(h)\r$ for $i\in K$. Here $v(h)\in L_K$ is 
regarded as lying in $L$. This observation 
leads to the following expression of $\hatvarv(\Delta,\V)$
which is sometimes useful.

\begin{equation}\label{eq:hatvarv2}
 \begin{split}
 \hatvarv(\Delta,\V) 
 = &\sum_{k=0}^n\sum_{K\in\sigmk,h_1\in \hat{H}_K}
   \zeta^{\l u^K, v(h_1)\r}
   \sum_{I\in \Sigma_K^{(n-k)}}\frac{w(I)}{|H_I|}\cdot \\
   &\sum_{h_2\in H_I}
    \prod_{i\in I\setminus K}
    \phi(-\l u_i^I, zv+v(h_2)\r, \tau, \sigma)
    \prod_{i\in K}
    \phi(-\l u_i^I, zv-\tau v(h_1)+v(h_2)\r, \tau, \sigma),
\end{split} 
\end{equation}
where $u^K=\sum_{i\in K}u_i^K$.

\begin{note}
In the sum above with respect to $K\in\sigmk$ and $h_1\in \hat{H}_K$, 
the term corresponding to $K=*\in \Sigma^{(0)}$ and $h_1=0\in \hat{H}_*=0$ 
is equal to $\varv(\Delta,\V)$.
\end{note}

It is sometimes useful as well to take a representative $v(h)$ of 
$h\in H_I$ such that
\begin{equation}\label{eqn:fIh}
 0\leq\l u_i^I,v(h)\r<1 \ \text{for all $i\in I$}. 
\end{equation}
Such a representative is unique. We denote the value 
$\l u_i^I,v(h)\r$ by $f_{I,h,i}$ for such a representative 
$v(h)$. If $h$ lies in $H_K$ for $K\in \sigmk$ contained in $I$, 
then $f_{I,h,i}=0$ for $i\not\in K$, and $f_{I,h,i}$ depends only 
on $K$ for $i\in K$ which we shall denote by $f_{K,h,i}$. 
The sum $\sum_{i\in K}f_{K,h,i}$ will be denoted by $f_{K,h}$. 
Note that \eqref{eqn:0} can be rewritten 
as
\begin{equation*}\label{eqn:0bis}
 \hat{H}_K=\{h\in H_K\mid f_{K,h,i}\not=0 
\quad \text{for any $i\in K$}. \}
\end{equation*} 
 
\begin{prop}\label{prop:laurent}
Let $\varv(\Delta,\V) =
\sum_{s=0}^\infty \var_{s}(z) q^s$ be the expansion
into power series, 
then $\zeta^{\frac{n}{2}}\var_{s}(z)$ belongs to 
$R(S^1)\otimes\Z[\zeta,\zeta^{-1}]$,
where $R(S^1)$ is identified with
$\Z[t,t^{-1}]$. Let $r$ be the least common multiple 
of $\{\vert H_I\vert\}_{I\in \sigmn}$. 
Then $\hatvarv(\Delta,\V)$ can be expanded in the form 
$\hatvarv(\Delta,\V) =
\sum_{s=0}^\infty \hatvar_{s}(z) q^s$, 
where $\zeta^{\frac{n}{2}}\hatvar_{s}(z)$ belongs to
$R(S^1)\otimes\Z[\zeta^{\frac1{r}},\zeta^{-\frac1{r}}]$.
\end{prop}
\begin{proof}
We introduce another variable $\tau_1$ with $\Im(\tau_1)>0$ and 
put $q_1=e^{2\pi\img \tau_1}$. Then, fixing $K\in \sigmk$ and 
$h_1\in \hatH_K$, we consider the function
\[
 \begin{split}
 \hatvar_{K,h_1}(z,&\tau_1,\tau,\sigma)= 
 \sum_{I\in \Sigma_K^{(n-k)}}\frac{w(I)}{|H_I|}\cdot \\
   &\sum_{h_2\in H_I}
    \prod_{i\in I\setminus K}
    \phi(-\l u_i^I, zv+v(h_2)\r, \tau, \sigma)
    \prod_{i\in K}
    \phi(-\l u_i^I, zv-\tau_1 v(h_1)+v(h_2)\r, \tau, \sigma).
 \end{split}
\]

Note that 
\begin{equation*}\label{eq:phiui2}
 \zeta^{\frac12}\phi(-\l u_i^I, zv-\tau v(h_1)+v(h_2)\r, \tau, \sigma)
 =\frac{1-\zeta{\xi_{i,h_2}^I}^{-1}}{1-{\xi_{i,h_2}^I}^{-1}}
\prod_{k=1}^{\infty}\frac{(1-\zeta{\xi_{i,h_2}^I}^{-1}q^k)
  (1-\zeta^{-1}\xi_{i,h_2}^Iq^k)}
 {(1-{\xi_{i,h_2}^I}^{-1}q^k)(1-\xi_{i,h_2}^Iq^k)}, 
\end{equation*}
where 
\[\begin{split} 
 \xi_{i,h_2}^I&=e^{2\pi\img\l u_i^I,zv-\tau_1v(h_1)+v(h_2)\r} \\
   &=\chi_I(u_i^I,h_2)t^{\l u_i^I,v\r}q_1^{-\l u_i^I,v(h_1)\r}. 
 \end{split} \] 
Since $\l u_i^I,v(h_1)\r=0$ for $i\in I\setminus K$, 
we have $\xi_{i,h_2}^I=\chi_I(u_i^I,h_2)t^{\l u_i^I,v\r}$ 
for $i\in I\setminus K$. 
We expand $\hatvar_{K,h_1}(z,\tau_1,\tau,\sigma)$ with respect 
to $\zeta$ and $q=e^{2\pi\img\tau}$ 
in the following form
\[ \hatvar_{K,h_1}(z,\tau_1,\tau,\sigma)= 
 \sum_{s_1\in \Z,s_2\in \Z_{\geq 0}}\hatvar_{s_1,s_2}(z,q_1)\zeta^{s_1}q^{s_2}. \]

Then we see that there is a family $\{\alpha_{I,j}\in \mathcal{O}_I\}_{j=1}^r$ 
for each $I$ such that $\tau_{I,I'}(\alpha_{I,j})=\alpha_{I',j}$ for 
$j=1,\ldots ,r$ and 
\[ \begin{split}
 \zeta^{\frac{n}{2}}\hatvar_{s_1,s_2}&(z,q_1)=\sum_{j=1}^r
  \sum_{I\supset K}\frac{w(I)}{|H_I|}\sum_{l\in \alpha_{I,j}}
  \sum_{h_2\in H_I}
 \frac{\prod_{i\in I}\chi_I(l,h_2)t^{\l l,v\r}q_1^{-\l l,v(h_1)\r}}
 {\prod_{i\in I}(1-{\xi_{i,h_2}^I}^{-1})} \\
 &=\sum_{j=1}^r
  \sum_{I\supset K}\frac{w(I)}{|H_I|}\sum_{l\in \alpha_{I,j}} 
  \sum_{h_2\in H_I}
 \frac{\prod_{i\in I}\chi_I(l,h_2)t^{\l l,v\r}q_1^{-\l l,v(h_1)\r}}
  {\prod_{i\in I\setminus K}(1-\chi_I(u_i^I,h_2)^{-1}t^{-\l u_i^I,v\r})
   \prod_{i\in K}(1-{\xi_{i,h_2}^I}^{-1})}. 
\end{split} \] 
We further expand each factor $\dfrac{1}{1-{\xi_{i,h_2}^I}^{-1}}$ 
in formal power series in the above expression. Then 
$\dfrac{1}{\prod_{i\in K}(1-{\xi_{i,h_2}^I}^{-1})}$ is expanded 
in an (infinite) sum of expressions of the form
$\sum_{l\in \beta_I}\chi_I(l,h_2)t^{\l l,v\r}q_1^{-\l l,v(h_1)\r}$, 
where $\beta_I$ is an element of $\mathcal{O}_I$ and it satisfies 
$\tau_{I,I'}(\beta_I)=\beta_{I'}$.
Hence $\zeta^{\frac{n}{2}}\hatvar_{s_1,s_2}(z,q_1)$ is written in a sum of 
expressions of the follwing form
\begin{equation}\label{eq:qone}
 \sum_{I\supset K}\frac{w(I)}{|H_I|}
 \sum_{l_1\in \alpha_{I,j},l_2\in \beta_I} 
  \sum_{h_2\in H_I}
 \frac{\prod_{i\in I}\chi_I(l_1+l_2,h_2)t^{\l l_1+l_2,v\r}
 q_1^{-\l l_1+l_2,v(h_1)\r}}
  {\prod_{i\in I\setminus K}(1-\chi_I(u_i^I,h_2)^{-1}t^{-\l u_i^I,v\r})}. 
\end{equation}

Since $\l u_i^I,v(h_1)\r=0$ for $i\in I\setminus K$, these 
expressions \eqref{eq:qone} 
belong to $\Z[t,t^{-1},q_1,q_1^{-1}]$ by Lemma \ref{lemm:integrality}. 
Hence $\zeta^{\frac{n}{2}}\hatvar_{s_1,s_2}(z,q_1)$ also belongs to 
$\Z[t,t^{-1},q_1,q_1^{-1}]$. Therefore, if we specialize 
$\tau_1$ to $\tau$ (hence $q_1$ to $q$) in 
$\hatvar_{K,h_1}(z,\tau_1,\tau,\sigma)= 
\sum_{s_1\in \Z,s_2\in \Z\geq 0}\hatvar_{s_1,s_2}(z,q_1)\zeta^{s_1}q^{s_2}$, 
we get an expansion 
\[ \hatvar_{K,h_1}(z,\tau,\tau,\sigma)= 
 \sum_{s\in \Z}\hatvar_{K,h_1,s}(z)q^{s} \]
with $\zeta^{\frac{n}{2}}\hatvar_{K,h_1,s}(z)\in 
\Z[t,t^{-1},\zeta,\zeta^{-1}]$. 
Then 
\[ \varv(\Delta,\V)= \sum_{K\in \Sigma,h_1\in \hatH_K}
 \zeta^{\l u^K,v(h_1)\r}\hatvar_{K,h_1}(z,\tau,\tau,\sigma)=
 \sum_{s\in \Z}\zeta^{\l u^K,v(h_1)\r}\hatvar_{K,h_1,s}(z)q^{s} \]
with $\zeta^{\frac{n}{2}}\zeta^{\l u^K,v(h_1)\r}\hatvar_{K,h_1,s}(z)\in 
\Z[t,t^{-1},\zeta^{\frac{1}{r}},\zeta^{-\frac{1}{r}}]$.
But $\varv(\Delta,\V)$ does not have negative powers 
of $q$ in its expansion as can be seen by 
taking a representative $v(h_1)$ of $h_1\in \hatH_I$ 
which satisfies \eqref{eqn:fIh} for each $I$ 
in the expression 
\eqref{eq:hatvarv} or \eqref{eq:hatvarv2}. 
It follows that 
$\varv(\Delta,\V)\zeta^{\frac{n}{2}}$ belongs to 
$(\Z[t,t^{-1},\zeta^{\frac{1}{r}},\zeta^{-\frac{1}{r}}] )[[q]]$.
This finishes the proof of Proposition \ref{prop:laurent}
for $\hatvarv(\Delta,\V)\zeta^{\frac{n}{2}}$. 

A similar easier argument shows that 
$\varv(\Delta,\V)\zeta^{\frac{n}{2}}$ 
belongs to $(R(S^1)\otimes \Z[\zeta,\zeta^{-1}])[[q]]$.
\end{proof}

The equivariant elliptic genus $\var(\Delta,\V)=\var(\Delta,\V;\tau,\sigma)
\in (R(T)\otimes \Z[\zeta,\zeta^{-1}])[[q]]$ 
and the equivariant orbifold elliptic genus $\hatvar(\Delta,\V)
=\hatvar(\Delta,\V;\tau,\sigma)\in 
(R(T)\otimes \Z[\zeta^{\frac1{r}},\zeta^{-\frac1{r}}])[[q]]$ are defined by
\[ v^*(\var(\Delta,\V))=\varv(\Delta,\V)\ \ \text{and}\ \ 
  v^*(\hatvar(\Delta,\V))=\hatvarv(\Delta,\V) \]
where one varies generic vectors $v$ in $L_{\V}$. 

\begin{rema}\label{rem:var}
Another but equivalent definition of $\var(\Delta,\V)$ is 
the following. Put $L_\C=L\otimes \C$ and let $T_\C=L_\C/L$ be 
the corresponding complex torus. 
The pairing $L^*\times L$ extends to a pairing $L^*\times L_\C$. 

Then 
\begin{equation*}\label{eq:var}
 \var(\Delta,\V)=\sum_{I\in\sigmn}\frac{w(I)}{|H_I|}
  \sum_{h\in H_I}
   \prod_{i\in I}\phi(\langle u_i^I,-w-v(h)\rangle,\tau,\sigma),
\end{equation*} 
where $w\in L_\C$. $\var(\Delta,\V)$ is a meromorphic function 
in $(w,\tau,\sigma)\in L_\C\times \mathcal{H}\times \C$. 
Similarly $\hatvar(\Delta,\V)$ may be defined by
\begin{equation*}\label{eq:hatvar}
 \begin{split}
 \hatvarv(\Delta,\V) 
 = &\sum_{I\in\sigmn}
   \frac{w(I)}{|H_I|}\cdot \\
   &\sum_{(h_1,h_2)\in H_I\times H_I}
  \prod_{i\in I}\zeta^{\l u_i^I, v(h_1)\r}
 \phi(\l u_i^I, -w+\tau v(h_1)-v(h_2)\r, \tau, \sigma).
\end{split} 
\end{equation*}
\end{rema}

The formulae in the following theorems express explicitly 
$\var(\Delta,\V)$ and $\hatvar(\Delta,\V)$ as
virtual characters of the torus $T$. When $\Delta$ is the fan
of a complete non-singular toric variety, the formula in 
Theorem \ref{thm:BL} is due to Borisov and Libgober \cite{BL1}.

\begin{theo}\label{thm:BL}
Let $\Delta=(\Sigma, C, w^\pm)$ be an $n$-dimensional complete simplicial 
multi-fan in a lattice $L$ and let $\V$ be as above. Then we have the equality
\[ \var(\Delta,\V)=\sum_{u\in L^*}t^{-u}\left(\sum_{k=0}^n\sum_{J\in \sigmk}
  (-1)^k\deg(\Delta_J)\prod_{j\in J}
  \frac1{1-\zeta q^{\langle u,v_j\rangle}}\right)\Phi(\sigma,\tau)^n ,\]
where 
$\deg(\Delta_J)$ is the degree of the projected multi-fan $\Delta_J$ of $J$
as defined in Section 2. 
\end{theo}
\begin{note}
The above equality is an equality of meromorphic functions as 
can be seen from the following proof. The right hand side 
is convergent to a meromorphic function in the domain
\[ \Im(\tau) >|\Im( \l u_i^I,w\r)|>0 \ \text{for all $i\in I$ 
and $I\in \sigmn$}. \]
Notice that $t^u$ is regarded as a function on $T_\C$ 
defined by $t^u(\exp(w))=e^{2\pi\img \l u,w\r}$, where 
the projeciton $L_\C\to T_\C$ is denoted by $\exp$ as usual. 
\end{note}
\begin{exam}
The equivariant elliptic genus of the complex projective space $\Proj^n$
is given by the following formula. The action of $T=T^n$ is the standard one 
given by 
\[ (g_1,\ldots ,g_n)[z_0,z_1,\ldots ,z_n]=[z_0,g_1z_1,\ldots ,g_nz_n]. \]
The lattice $L$ is identified with $\Z^n$ and 
$\V=\{v_i\}$ is given by 
\[ v_i=e_i\ \text{(standard unit vector) for $i=1,\ldots ,n$\ and 
$v_{n+1}=-(e_1+\cdots +e_n)$}.\] 
\[
\var(\Proj^n)=\left(\sum_{u=(m_1,\ldots,m_n)\in\Z^n}t^{-u}
 \frac{1-\zeta^{n+1}}{(1-\zeta q^{m_1})\cdots
 (1-\zeta q^{m_n})(1-\zeta q^{-m_1\cdots -m_n})}\right)
 (-\Phi(\sigma,\tau))^n,
\]
where $L^*=H^2(BT)$ is identified with $\Z^n$.
\end{exam}

\begin{theo}\label{thm:BL2}
Let $(\Delta,\V)$ be as in Theorem \ref{thm:BL}. Then we have 
the equality 
\begin{align*}
 &\hatvar(\Delta,\V) \\ 
 =&\sum_{u\in L^*}t^{-u}\left(\sum_{k=0}^n
  \sum_{J\in\Sigma^{(k)}, h\in H_J}(-1)^k  
  \deg(\Delta_J)\zeta^{f_{J,h}}q^{\l u,v_{J,h}\r}
  \prod_{i\in J}\frac1{1-\zeta q^{\l u,v_i \r}}\right)\Phi(\sigma,\tau)^n,
\end{align*}
where $v_{J,h}=\sum_{i\in J}f_{J,h,i}v_i$.
\end{theo}
\begin{note}
It can be shown that $v_{J,h}$ belongs to $L_J$ and 
hence $\l u,v_{J,h}\r\in \Z$. This fact shows that 
$\hatvar(\Delta,\V)$ belongs to
$(R(T)\otimes \Z[\zeta^{\frac1{r}},\zeta^{-\frac1{r}}])\Z[[q]]$. This fact 
was aleady proved in Proposition \ref{prop:laurent}. 
We also see that 
\[\zeta^{f_{J,h}}q^{\l u,v_{J,h}\r}
  \prod_{i\in J}\frac1{1-\zeta q^{\l u,v_i \r}}=\prod_{i\in J}
  \frac{(\zeta q^{\l u,v_i\r})^{f_{J,h,i}}}{1-\zeta q^{\l u,v_i \r}}. \] 
\end{note}

\begin{exam}
Let $b>1$ be an integer. Let $H=\{g\in S^1\vert g^b=1\}$ act on 
$\Proj^2$ by $g[z_0,z_1,z_2]=[z_0,z_1,gz_2]$   
and put $M=\Proj^2/H$. The action of $T=T^2$ on $M$ is induced from 
the action on $\Proj^2$ given by
\[ (g_1,g_2)[z_0,z_1,z_2]=[z_0,g_1z_1,g_2z_2]. \]
The lattice $L$ is identified with $\Z^2$ and $\V=\{v_i\}$ is given by 
\[ v_1=e_1,v_2=be_2, v_3=-(e_1+be_2). \]
Then a calculation using Theorem \ref{thm:BL2} 
yields
\[ \hatvar(M)=\sum_{(m_1,m_2)\in\Z^2}a_{m_1,m_2}t_1^{-m_1}t_2^{-m_2} \]
with
\[ a_{m_1,m_2}=\frac{(1-\zeta^{3/b})(1-\zeta^2q^{-bm_2})\Phi(\sigma,\tau)^2}
 {(1-\zeta q^{m_1})(1-\zeta q^{-m_1-bm_2})
 (1-\zeta^{1/b}q^{m_2})(1-\zeta^{2/b}q^{-m_2})}.
\]
\end{exam}

For the proof of Theorems \ref{thm:BL} and \ref{thm:BL2} 
we need lemmas below. 

\begin{lemm}\label{lemm:BL1}
Suppose that $|q|<|t|<1$. Then we have the equality
\begin{equation}\label{eq:BL1}
 \phi(z,\tau,\sigma)=-\Phi(\sigma,\tau)
\sum_{m\in \Z}\frac{t^m}{1-\zeta q^m}.
\end{equation}
\end{lemm}

\begin{lemm}\label{lemm:BL}
Put $\alpha=e^{2\pi\img w}$. Suppose that $|\alpha|=1$ 
and $|q|<|t|<1$. 
If $l\not=0$ is an integer, then we have the equality
\[ \phi(lz+w,\tau,\sigma)=
\begin{cases}
 -\Phi(\sigma,\tau)\sum_{m\in \Z}\alpha^mt^{lm}\dfrac1{1-\zeta q^m} & 
 \text{if \ $l>0$}, \\
 -\Phi(\sigma,\tau)\sum_{m\in \Z}\alpha^mt^{lm}
 \left(\dfrac1{1-\zeta q^m}-1\right) & \text{if \ $l<0$}. 
\end{cases} \]
\end{lemm}

\begin{proof}
First we prove \eqref{eq:BL1}. The argument follows that of \cite{BL1}. 
We repeat their argument for the sake of completeness. 
Recall that $\zeta=e^{2\pi\img\sigma}$. We regard the both sides of 
\eqref{eq:BL1} as meromorphic functions of $\sigma$ and consider 
the quotient 
\[
  f(\sigma)=\phi(z,\tau,\sigma)\slash\left(-\Phi(\sigma,\tau)
   \sum_{m\in \Z}\dfrac{t^m}{1-\zeta q^m}\right) .
\]
defined for $|q|<|t|<1$. It is doubly periodic in $\sigma$ 
with respect to the lattice $\Z\tau +\Z$ and has no poles.
Its value at $\sigma=0$ is $1$. This implies the equality \eqref{eq:BL1}.

We now prove Lemma \ref{lemm:BL}. So suppose $|\alpha|=1$ and
$|q|<|t^{|l|}|<1$. If $l>0$ then $|q|<|\alpha t^l|<1$ and the equality in
Lemma \ref{lemm:BL} is \eqref{eq:BL1} itself since 
$e^{2\pi\img(lz+w)}=\alpha t^l$.
Next we consider the case $l<0$. Since  $|\alpha t^l|>1$, 
Lemma \ref{lemm:BL1} can not be applied.
So we proceed in the following way. First observe that
\[ \phi(-z,\tau,\sigma)=\phi(z,\tau,-\sigma)\ \ \text{and}\ \ 
  \Phi(\sigma,\tau)=-\Phi(-\sigma,\tau). \] 
Hence, if we put $l'=-l$, we see easily that
\[\phi(lz+w,\tau,\sigma)/\Phi(\sigma,\tau)=
 -\phi(l'z-w,\tau,-\sigma)/\Phi(-\sigma,\tau). \]
We then apply Lemma \ref{lemm:BL1} to the right hand side and get
\begin{align*}
 \phi(lz+w,\tau,\sigma)/\Phi(\sigma,\tau)
  &= \sum_{m\in \Z}\alpha^{-m}t^{l'm}
 \dfrac1{1-\zeta^{-1} q^m} \\
 &=-\sum_{m\in \Z}\alpha^mt^{lm}
 \left(\dfrac1{1-\zeta q^m}-1\right). 
\end{align*}
This proves Lemma \ref{lemm:BL}.
\end{proof}

We now proceed to the proof of Theorem \ref{thm:BL}. Take a generic vector
$v\in L_\V$. Then $\langle u_i^I,v\rangle$ is an integer
for any $I\in\sigmn$ and $i\in I$.  
For $I\in \sigmn$ we put
$I(v)=\{i\in I|\langle u_i^I,v\rangle <0\}$. 
Suppose that $|q|<|t^{|\l u_i^I,v\r|}|<1$ for all $I\in \sigmn$ and $i\in I$ .
Using Lemma \ref{lemm:BL}
we have
\begin{multline*}
 \frac1{(-\Phi(\sigma,\tau))^n}\prod_{i\in I}
\phi(-\l u_i^I,zv+v(h)\r,\tau,\sigma)  \\
\shoveleft{    
=\prod_{i\in I(v)}\left(\sum_{m_i\in \Z}\chi_I(u_i^I,h)^{-m_i}
 t^{-m_i\langle u_i^I,v\rangle}
 \frac1{1-\zeta q^{m_i}}\right) } \\
\prod_{i\in I\setminus I(v)}\left(\sum_{m_i\in \Z}\chi_I(u_i^I,h)^{-m_i}
 t^{-m_i\langle u_i^I,v\rangle}
\left(\frac1{1-\zeta q^{m_i}}-1\right)\right) \\   
\shoveleft{
=\sum_{m_i\in\Z, i\in I}
\left(\prod_{i\in I}\chi_I(u_i^I,h)^{-m_i}
    \prod_{i\in I}t^{-m_i\langle u_i^I,v\rangle}
    \prod_{i\in I(v)}\frac1{1-\zeta q^{m_i}}\prod_{i\in {I\setminus I(v)}}
\left(\frac1{1-\zeta q^{m_i}}-1\right)
\right)}
\end{multline*}

If we put $u=\sum_{i\in I}m_iu_i^I \in L_{I,\V}^*$, 
then $m_i=\langle u,v_i\rangle$.
Therefore $\prod_{i\in I}t^{-m_i\langle u_i^I,v\rangle}
=t^{-\langle u,v\rangle}$.
Since $\chi_I(u,\ )=e^{2\pi\sqrt{-1}\langle u,\ \rangle}$ we see that 
$\prod_{i\in I}\chi_I(u_i^I,h)^{-m_i}=\chi_I(u,h)^{-1}$. Since 
the $1$-dimensional representation of $H_I=L/L_{I,\V}$ is trivial if 
and only if $u\in L^*$, it follows that
\[ \sum_{h\in H_I}\chi_I(u,h)^{-1}= 
\begin{cases}
  |H_I|\quad &\text{if}\ u\in L^*, \\
  0 &\text{if}\ u\not\in L^*. 
\end{cases} \]
Furthermore we see easily that
\begin{equation*}\label{eq:sumJ}
\prod_{i\in I(v)}\frac1{1-\zeta q^{m_i}}
\prod_{i\in {I\setminus I(v)}}\left(\frac1{1-\zeta q^{m_i}}-1\right)
=\sum_{k=0}^n\sum_{J\in \sigmk :I(v)\subset J\subset I}(-1)^{n-k}
\prod_{j\in J}\frac1{1-\zeta q^{m_j}},
\end{equation*}
for each $I\in \sigmn$. Combining these we have
\begin{multline*}
 \sum_{h\in H_I}\frac1{\Phi(\sigma,\tau)^n}
 \prod_{i\in I}\phi(-\l u_i^I,zv+v(h)\r,\tau,\sigma) \\   
 =|H_I|\sum_{u\in L^*}t^{-\langle u,v\rangle}\left(
   \sum_{k=0}^n\sum_{J\in \sigmk :I(v)\subset J\subset I}(-1)^k
   \prod_{j\in J}\frac1{1-\zeta q^{\langle u,v_j\rangle}}\right).
\end{multline*}

Finally, applying this to (\ref{eq:varv}) we have
\begin{multline*}
\frac1{\Phi(\sigma,\tau)^n}v^*\left(\var(\Delta,\V)\right) 
=\frac1{\Phi(\sigma,\tau)^n}\varv(\Delta,\V) \\
\shoveleft{
=\sum_{u\in L^*}t^{-\langle u,v\rangle}\left(
   \sum_{I\in \sigmn}w(I)\left(
  \sum_{k=0}^n\sum_{J\in \sigmk :I(v)\subset J\subset I}(-1)^k
   \prod_{j\in J}\frac1{1-\zeta q^{\langle u,v_j\rangle}}
  \right)\right)} \\
\shoveleft{
=\sum_{u\in L^*}t^{-\langle u,v\rangle}\left(
\sum_{k=0}^n\sum_{J\in \sigmk}(-1)^k
\left(\sum_{I\in \sigmn :I(v)\subset J\subset I}w(I)\right)
   \prod_{j\in J}\frac1{1-\zeta q^{\langle u,v_j\rangle}}\right)
}\\
\shoveleft{
=\sum_{u\in L^*}t^{-\langle u,v\rangle}\left(
\sum_{k=0}^n\sum_{J\in \sigmk}(-1)^k
 \deg(\Delta_J)
 \prod_{j\in J}\frac1{1-\zeta q^{\langle u,v_j\rangle}}\right),
} \\
\end{multline*}
since $\sum_{I\in \sigmn :I(v)\subset J\subset I}w(I)=
\deg(\Delta_J)$ by definition.
Since $\var(\Delta,\V)$ belongs to $(R(T)\otimes\Z[\zeta,\zeta^{-1}])[[q]]$
and the last equality holds for any generic vector $v$,  
Theorem \ref{thm:BL} follows. 

We next prove Theorem \ref{thm:BL2}. We use the following
\begin{lemm}\label{lemm:1}
Let $f$ be a real number with $0<f<1$ and $l\not=0$ an integer.
If
\[ |q^f|,|q^{1-f}|<|t^{|l|}|,\ |t|\leq 1 \]
then 
\begin{equation*}\label{eqn:1}
 \phi(lz+f\tau+w,\tau,\sigma)=
 -\Phi(\sigma,\tau)\sum_{m\in \Z}\alpha^m t^{lm}
 \frac{q^{fm}}{1-\zeta q^m},
\end{equation*}
where $\alpha=e^{2\pi\img w},\ |\alpha|=1$ as before.
\end{lemm}
\begin {proof}
In view of Lemma \ref{lemm:BL1} it is enough to show that
\[ |q|<|q^ft^l|<1 \]
regardless of the sign of $l$.

Suppose that $l>0$. Then $|q^{1-f}|<|t^l|$ implies $|q|<|q^ft^l|$.
Since $|q^f|<1$ and $|t^l|\leq 1$, we have $|q^ft^l|<1$.

Soppose that $l<0$. Then $|q^f|<|t^{-l}|$ implies $|q^ft^l|<1$. 
Also $|q^{1-f}|<|t^{-l}|$ implies $|q|<|q^ft^{-l}|$. But 
$|q^ft^{-l}|\leq |q^ft^l|$ since $l<0$ and $|t|\leq 1$. 
Hence $|q|<|q^ft^l|$.  
\end{proof}

We now proceed to the proof of Theorem \ref{thm:BL2}. Take a generic vector
$v\in L_\V$.  
Fix  $i\in K$ and take the representative $v(h_1)$ of $h_1\in \hat{H}_K$ 
such that 
\[ \l u_i^I,v(h_1)\r=\l u_i^K,v(h_1)\r =f_{K,h_1,i}.\]
If $t\in \C$ satisfies
\[ |q^{f_{K,h_1,i}}|,|q^{1-f_{K,h_1,i}}|<|t^{|\uv|}|, \quad |t|<1, \]
then by Lemma \ref{lemm:1}, we have
\[ \begin{split}
 \frac1{-\Phi(\sigma,\tau)}&\phi(-\l u_i^I,zv-\tau v(h_1)+v(h_2)\r,
\tau, \sigma) \\
&=\sum_{m_i\in\Z}\chi_I(u_i^I,v(h_2))^{-m_i}q^{m_if_{K,h_1,i}}t^{-m_i\uv}
\frac1{1-\zeta q^{m_i}}. 
\end{split} \]

Next fix $i\in I\setminus K$. Then, by Lemma \ref{lemm:BL}, we have
\[ \begin{split}
 &\frac1{-\Phi(\sigma,\tau)}\phi(-\l u_i^I,zv+v(h_2)\r,\tau,\sigma) \\ 
 &= \begin{cases}
  \sum_{m_i\in\Z}\chi_I(u_i^I,v)^{-m_i}t^{-m_i\uv}\dfrac1{1-\zeta q^{m_i}}
  \qquad \text{for} 
  \ i\in I(v)\setminus K, \\
  \sum_{m_i\in\Z}\chi_I(u_i^I,v)^{-m_i}t^{-m_i\uv}
 \left(\dfrac1{1-\zeta q^{m_i}}-1\right) 
  \quad \text{for} \ i\in I\setminus (I(v)\cup K).
\end{cases}
\end{split}
\]

Now suppose that $|t|<1$ and $q$ satisfy 
\[ |q^{f_{K,h_1,i}}|, |q^{1-f_{K,h_1,i}}|<|t^{|\l u_i^I,v\r|}| \]
for all $i\in K,\ h_1\in \hat{H}_K, \ K\in \sigmk$, 
and
\[ |q|<|t^{|\l u_i^I,v\r|}| \]
for all $i\in I\setminus K, I\in \Sigma_K^{(n-k)}$ and 
$K\in \sigmk$. Then we obtain

\begin{equation*}
\begin{split}
&\frac{1}{(-\Phi(\sigma,\tau))^n}\prod_{i\in I\setminus K}
\phi(-\l u_i^I,zv+v(h_2)\r,\tau,\sigma)
  \prod_{i\in K}\phi(-\l u_i^I,zv-\tau v(h_1)+v(h_2)\r,\tau,\sigma) \\
&=\prod_{i\in I(v)\setminus K}{\sum_{m_i\in \Z}
\chi_I(u_i^I,h_2)^{-m_i}t^{-m_i\uv}\frac{1}{1-\zeta q^{m_i}}} \\
&\quad \prod_{i\in I\setminus (I(v)\cup K)}{\sum_{m_i\in \Z}
 \chi_I(u_i^I,h_2)^{-m_i}
 t^{-m_i\uv}\left(\frac{1}{1-\zeta q^{m_i}}-1\right)} \cdot \\
&\qquad \prod_{i\in K}{\sum_{m_i\in \Z}
\chi_I(u_i^I,h_2)^{-m_i}q^{m_if_{K,h_1,i}}
 t^{-m_i\uv}\frac{1}{1-\zeta q^{m_i}}}.
\end{split}
\end{equation*}

Using the above equality and \eqref{eq:hatvarv2}, and 
arguing as in the proof of Theorem \ref{thm:BL} 
we have
\begin{equation}\label{eqn:2}
 \begin{split}
 &\hatvar(\Delta,\V)/\Phi(\sigma,\tau)^n \\ 
 =&\sum_{u\in L^*}t^{-u}\left(\sum_{k=0}^n\sum_{K\in\Sigma^{(k)},h\in \hat H_K}
 \zeta^{f_{K,h}}q^{\l u,v_{K,h}\r}\sum_{K\subset J}(-1)^k  
  \deg(\Delta_J)
  \prod_{j\in J}\frac{1}{1-\zeta q^{\l u,v_j \r}}\right).
 \end{split}
\end{equation}
Fix $J\in \Sigma$. It is easy to see that the union of 
$\{\hat{H}_K\mid K\in \Sigma, K\subset J\}$ is disjoint. 
Since any $h\in H_J$ is contained in 
$\hat{H}_{K_h}$ by (\ref{eqn:0}) where 
$K_h=\{j\in J\mid f_{J,h,j}\not=0\}$, we have
$H_J=\sqcup \hat{H}_K$.
Moreover we have
\[ f_{J,h}=f_{K_h,h} \quad\text{and}\quad v_{J,h}=v_{K_h,h}.\]
Taking these facts in accout the right hand side of the equality
(\ref{eqn:2}) is transformed into
\[ \sum_{u\in L^*}t^{-u}\left(\sum_{k=0}^n
\sum_{J\in\Sigma^{(k)},h\in H_J}(-1)^k  
  \deg(\Delta_J)\zeta^{f_{J,h}}q^{\l u,v_{J,h}\r}
  \prod_{i\in J}\frac1{1-\zeta q^{\l u,v_i \r}}\right). \]
This proves Theorem \ref{thm:BL2}.

The elliptic genus $\var(\Delta,\V)$ reduces to the 
so-called $T_y$-genus for $q=0$ if it is multiplied by $\zeta^{n/2}$ 
and if $\zeta$ is substituted by $-y$. 
Namely 
\[ T_y(\Delta,\V)=\sum_{I\in \sigmn}\frac{w(I)}{|H_I|}
   \sum_{h\in H_I}\prod_{i\in I}
   \frac{1+y\chi_I(u_i^I,h)^{-1}t^{-u_i^I}}
 {1-\chi_I(u_i^I,h)^{-1}t^{-u_i^I}}. \]
In \cite{HM} it was shown that the equivariant $T_y$-genus of
a torus manifold was rigid. The same proof is valid for general 
complete simplicial multi-fans. For the sake of 
completeness we review the argument briefly. 
Let $v$ be a generic vector. We set
\[ \mu(I)=\#\{i\in I\mid \l u_i^I,v\r >0,\} \]
for $I\in \sigmn$ and define 
\[ h_k(\Delta):=\sum_{I\in \sigmn,\ \mu(I)=k}w(I), \]
for each integer $k$ with $0\leq k\leq n$. 

We consider $v^*(T_y(\Delta,\V))$. It is written in the form 
$\sum_{I\in \sigmn}R_I(t)$, where
\[ R_I(t)=\frac{w(I)}{|H_I|}
   \sum_{h\in H_I}\prod_{i\in I}
   \frac{1+y\chi_I(u_i^I,h)^{-1}t^{-\l u_i^I,v\r}}
        {1-\chi_I(u_i^I,h)^{-1}t^{-\l u_i^I,v\r}}. \]
Regarded as a rational function of $t$, $R_I(t)$ takes value 
\[ w(I)(-y)^{\mu(I)}, \]
at $t=0$ and similarly
\[ w(I)(-y)^{n-\mu(I)} \]
at $t=\infty$. Hence $v^*(T_y(\Delta,\V))(t)$ takes finite values 
at $t=0$ and $t=\infty$. On the other hand it is 
a Laurent polynomial in $t$. Hence it must be a constant. 
Since this is true for any generic vector $v$, $T_y(\Delta,\V)$ is
a constant. Moreover that constant is equal either to 
\[ \sum_{I\in \sigmn}w(I)(-y)^{\mu(I)}=\sum_{k=0}^nh_k(\Delta)(-y)^k \]
or to 
\[ \sum_{I\in \sigmn}w(I)(-y)^{n-\mu(I)}=\sum_{k=0}^nh_{n-k}(\Delta)(-y)^k. \] 
It follows that $h_k(\Delta)$ does not depend on $v$ 
and that $h_k(\Delta)=h_{n-k}(\Delta)$. 
Moreover we see that $T_y(\Delta,\V)$ is independent of $\V$, 
since $h_k(\Delta)$ 
depends only on $\Delta$.  So we call it the $T_y$-genus of $\Delta$ and 
simply write $T_y[\Delta]$. Thus we have proved
\begin{prop}\label{prop:Tygenus}
\begin{equation}\label{eq:h}
 T_y[\Delta]=\sum_{k=0}^nh_k(\Delta)(-y)^k . 
\end{equation}
Here the equality $h_{n-k}(\Delta)=h_k(\Delta)$ holds.
\end{prop}
\begin{note}
$h_0(\Delta)=T_0[\Delta]$ is the Todd genus of $\Delta$, 
and $h_n(\Delta)=\deg(\Delta)$ by definition of the latter. 
Hence $T_0[\Delta]$ equals $\deg(\Delta)$, cf. \cite{HM}.
\end{note}

The following Proposition was proved in \cite{HM}. 
We shall give a different proof using Theorem \ref{thm:BL}. 
\begin{prop}\label{prop:T_y}
\begin{equation}\label{eq:T_y}
 T_y[\Delta]=\sum_{k=0}^ne_k(\Delta)(-1-y)^{n-k}
\end{equation}
where $e_k(\Delta)=\sum_{J\in \sigmk}\deg(\Delta_J)$.
\end{prop}
\begin{proof} 
We look at the coefficient of $t^u$ for $u=0$ and with $q=0$ 
in Theorem \ref{thm:BL} which is equal to $T_y[\Delta]\zeta^{-\frac{n}{2}}$. 
Noting that $\Phi(\sigma,\tau)$ is approaching to $(-1-y)\zeta^{-\frac12}$ 
when $\Im \tau$ approaches to $\infty$, we
obtain \eqref{eq:T_y}.
\end{proof}


\section{Equivariant first Chern class}\label{sec:chern}
Let $\Delta$ be a complete simplicial multi-fan 
in a lattice $L$ and $\V=\{v_i\}_{i\in \sigmone}$ a set of prescibed 
edge vectors as before. An $H^*(BT)$-module structure of 
$H_T^*(\Delta)$ is defined by \eqref{eq:structure}. The class 
\[ \sum_{i\in \sigmone}x_i \in H_T^2(\Delta) \]
will be called the \emph{equivariant first Chern class} of the 
pair $(\Delta,\V)$, and 
will be denoted by $c_1^T(\Delta,\V)$. When $\Delta$ is non-singular,
$\V$ consists of primitive vectors which is determined by $\Delta$ by
our convention. In this case we simply write $c_1^T(\Delta)$ 
and call it the \emph{equivariant first Chern class} of the non-singular multi-fan 
$\Delta$.

The image of $c_1^T(\Delta,\V)$ in
$H_T^2(\Delta)/H^2(BT)$ is called the first Chern class
of $(\Delta,\V)$ and is denoted by $c_1(\Delta,\V)$.
Let $N>1$ be an integer. The first Chern class $c_1(\Delta,\V)$ is
divisible by $N$ if and only if $c_1^T(\Delta,\V)$ is of the form
\[ c_1^T(\Delta,\V)=Nx +u,\ x\in H_T^2(\Delta),\ u\in H^2(BT). \]
We set $u^I=\iota_I^*(c_1^T(\Delta,\V))=\sum_{i\in I}u_i^I\in L_{I,\V}^*$.
Note that $u^I$ does not belong to $L^*=H^2(BT)$ in general.

\begin{lemm}\label{lemm:Ndivisible}
The following three conditions are equivalent:
\begin{enumerate}
\item the first Chern class $c_1(\Delta,\V)$ is
divisible by $N$, 
\item $u^I \bmod N $ regarded as an element of 
$L_{\V}^*/NL_{\V}^*$ is independent of $I\in \sigmn$ and 
belongs to the image of $L^*=H^2(BT)$,
\item there is an element $u\in H^2(BT)$ such that 
$\langle u,v_i\rangle =1\bmod N$ for all $i\in \sigmone$.
\end{enumerate}
\end{lemm}
\begin{proof}
Suppose $c_1^T(\Delta,\V)$ is of the form
\[ c_1^T(\Delta,\V)=Nx +u,\ x\in H_T^2(\Delta),\ u\in H^2(BT). \]
Then $u^I=\iota_I^*(c_1^T(\Delta,\V)) \bmod N$ is equal to 
$\iota_I^*(u)=u \bmod NL_{\V}^*$, and hence belongs to 
the image of $H^2(BT)$. Thus (i) implies (ii). 

Suppose that $u\in H^2(BT)$ and $u\bmod N $ is equal to $u^I\bmod N$ for any 
$I\in \sigmn$. Then
$\langle u,v_i\rangle=\langle u^I,v_i\rangle \bmod N$, and 
hence 
$\langle u,v_i\rangle=\sum_{j\in I}\langle u_j^I,v_i\rangle=1 \bmod N$
for $i\in I$. Thus (ii) implies (iii).

Suppose $\langle u,v_i\rangle=1 \bmod N$ for any $i\in \sigmone$. 
Then, by (\ref{eq:structure}),
\[ c_1^T(\Delta,\V)-u=\sum_{i\in \sigmone}(1-\langle u,v_i\rangle)x_i
  =0 \bmod N H_T^2(\Delta) .\]
Hence $c_1^T(\Delta,\V)$ is of the form $c_1^T(\Delta,\V)=Nx+u$ with 
$u\in H^2(BT)$. Thus (iii) implies (i).
\end{proof}

\begin{rema}\label{rem:face ring}
Let $M$ be a torus manifold
and $\Delta(M)$ its associated multi-fan. 
Put $\hat{H}_T^2(M)=H_T^2(M)/S\text{-torsion}$. In \cite{M} 
it was shown that there is a canonical
embedding of $H_T^2(\Delta(M))$ in $\hat{H}_T^2(M)$, and, in case
$M$ is a stably almost complex torus manifold, 
$c_1^T(M)\in H_T^2(M)$ descends to $c_1^T(\Delta(M))\in H_T^2(\Delta(M))$. 
It follows that,
if $M$ is a stably almost complex torus manifold and $c_1(M)$ is divisible by
$N$, then $c_1(\Delta(M))$ is also divisible by $N$. Even if $M$ is
a stably almost complex orbifold it can be shown that 
$c_1^T(M)\in H_T^2(M;\R)$ descends to 
$c_1^T(\Delta(M),\V(M))\in H_T^2(\Delta(M))\otimes\R$. But the divisibility of
the first Chern class has no meaning with real coefficients. 
We have to work with orbifold cohomology theory with integer coefficents. 
We understand it in this paper
by that of $c_1(\Delta(M),\V(M))$. 
\end{rema}

The following property (P) will be called the global type condition:
\vspace*{0.2cm}\\
(P):\  $L_{I,\V}=L_\V$ for all $I\in \sigmn$. 
\vspace*{0.2cm}\\
Typical examples of pairs $(\Delta,\V)$ satisfying the condition (P) 
are provided by global torus
orbifolds. In fact, let $\tilde{M}$ be a torus manifold of diemnsion 
$2n$ and let $H$ be
a finite subgroup of the torus $\tilde{T}$ acting on $\tilde{M}$.
The quotient $M=\tilde{M}/H$ is a torus orbifold equipped with 
an orbifold structure for which $(\tilde{M},M,H,\pi)$ is an orbifold 
chart where $\pi:\tilde{M}\to M$ is the projection. If $\tilde{M}_i$
is a characteristic submanifold of $\tilde{M}$, then its image $M_i$ 
by $\pi$ is a characteristic suborbifold of $M$. Conversely every 
characteristic suborbifold of $M$ is of the above form. It follows that one 
can identify the simplicial set $\Sigma(M)$ with $\Sigma(\tilde{M})$.
The lattice $\tilde{L}$ for the multi-fan $\Delta(\tilde{M})$ 
is identified with $\pi_1(\tilde{T})$ and similarly the multi-fan
$\Delta(M)$ is defined in the lattice $L=\pi_1(T)$. 
Let $\{\tilde{v}_i\}_{i\in \sigmone(\tilde{M})}$ be the primitive generators
corresponding to the oriented characteristic submanifolds 
$\tilde{M}_i$. Put $v_i=\pi_*(\tilde{v}_i)\in L$ and
$\V=\{v_i\}$. Then $L_\V$ coincides with the image $\pi_*(\tilde{L})$. 
For any $I\in \sigmn(\tilde{M})$ the lattice $\tilde{L}$ 
is generated by $\{\tilde{v}_i\mid i\in I\}$. Hence $L_\V$ is generated 
by $\{v_i\mid i\in I\}$. This shows that the pair $(\Delta(M),\V)$ satisfies 
the condition (P).
Note that, when $\Delta$ is a non-singular multi-fan, we have
$L_{I,\V}=L_\V=L$ for all $I\in \sigmn$ and the condition (P) is
automatically satisfied.

\begin{rema}\label{rem:cover}
If $(\Delta,\V)$ satisfies the condition (P), then there is a complete 
non-singular multi-fan $\tilde{\Delta}=
(\tilde{\Sigma},\tilde{C},\tilde{w}^\pm)$ 
in a lattice $\tilde{L}$ and an injective homomorphism 
$\pi: \tilde{L}\to L$ such that $\pi(\tilde{L})=L_\V$. Namely let 
$\tilde{L}$ be a copy of $L_\V$ and let $\tilde{v}_i$ be the copy of $v_i$ 
in $\tilde{L}$. The identification of $\tilde{L}$ with $L_\V$ composed with 
the inclusion of $L_\V$ in $L$ defines the map $\pi$. We put 
$\tilde{\Sigma}=\Sigma,\ \tilde{w}^\pm=w^\pm$ and define $\tilde{C}(I)$ 
to be the cone generated by $\{\tilde{v}_i\}_{i\in I}$. These define
the multi-fan $\tilde{\Delta}$. We may call it the (ramified)
covering of $\Delta$ with respect to $\mathcal{V}$. 
Since $\{\tilde{v}_i\mid i\in I\}$ is a basis of $\tilde{L}$ for any 
$I \in \Sigma^{(n)}$, $\tilde{\Delta}$ is non-singular. 
We put $\tilde{T}=\tilde{L}_\R/\tilde{L}$. It is the torus associated
with the multi-fan $\tilde{\Delta}$. Then the kernel of the induced
map $\pi:\tilde{T} \to T=L_\R/L$ is identified with $H=L/L_\V$.
Conversely given a complete non-singular multi-fan 
$\tilde{\Delta}=(\tilde{\Sigma},\tilde{C},\tilde{w}^\pm)$ 
in a lattice $\tilde{L}$ and 
a sublattice $L$ of $\tilde{L}_\R$ such that
$\tilde{L}\subset L$ we can construct
a new complete simplicial multi-fan $\Delta$ 
in the lattice $L$ 
and a set of edge vectors $\mathcal{V}$ such that the covering of 
$\Delta$ with respect to $\mathcal{V}$ coincides with $\tilde{\Delta}$.
Geometric picture of this construction is making quotient torus
orbifold $\tilde{M}/H$ out of a torus manifold $\tilde{M}$ by 
a subgroup $H$ of 
the torus $\tilde{T}$ acting on $\tilde{M}$.
\end{rema}
\begin{rema}\label{rem:N-divisible}
The following fact can be proved easily. 
Under the situation of Remark \ref{rem:cover} suppose that 
$c_1(\Delta, \V)$ is divisible by $N$, then $c_1(\tilde{\Delta})$ 
is also divisible by $N$. Conversely if $c_1(\tilde{\Delta})$ is
divisible by $N$ and the order of $H=L/L_\V$ is 
relatively prime to $N$, then 
$c_1(\Delta)$ is also divisible by $N$.
\end{rema}

Let $v$ be a generic vector in $L_\V\subset L=H_2(BT)$ . If we 
fix $I$ and write 
\[ v=\sum_{i\in I}m_iv_i, \]
then $m_i=\langle u_i^I,v\rangle\in \Z $. Let $m\geq 1$ be an 
integer. Fixing $v$ and $m$, we put
\[ I_{(m)}:=\{i\in I\mid m \ \text{does not divide}\ m_i \} \]
for $I\in \sigmn$. It will be called $\bmod\: m$ face of $I$. 
Note that $I_{(m)}$ depends on $v$.

\begin{lemm}\label{lemm:modmface}
Suppose that the condition (P) is satisfied. Let $v,m$ and $I$ be 
as above, and let $K=I_{(m)}$ be the $\bmod\: m$
face of $I$. If $I'\in \sigmn$ contains $K$, then $K$ is also
the $\bmod\: m$ face of $I'$. Moreover, $\langle u_i^{I'},v\rangle
= \langle u_i^{I},v\rangle \bmod\: m$ for $i\in K$.
\end{lemm}
\begin{proof}
We put $m_i=\langle u_i^I,v\rangle$ and 
$m'_{i'}=\langle u_{i'}^{I'},v\rangle$. Then we have
\[ v=\sum_{i\in I}m_iv_i=\sum_{i'\in I'}m'_{i'}v_{i'}.\]
By assumption $m_i=0 \bmod\: m$ for $i\notin K$. Hence
\[ \sum_{i\in K}m_iv_i=\sum_{i'\in I'}m'_{i'}v_{i'}\ \ \bmod\: m \ 
\text{in $L_\V$}\]
or 
\[ \sum_{i\in K}(m'_i-m_i)v_i+
    \sum_{i'\in I',\ i'\notin K}m'_{i'}v_{i'}=0 \ \ \bmod\: m \ 
\text{in $L_\V$}.\]
Since $\{v_{i'}\mid i'\in I'\}$
is a basis of the free module $L_{I',\V}=L_\V$, we see that
$m'_{i'}=0 \bmod\: m$ for $i'\in I',\ i'\notin K$ and $m'_i=m_i \bmod\: m$ 
for $i\in K$.
\end{proof}
We shall say that $I$ and $I'$ are $(v,m)$-equivalent and
write $I\sim I'$ if 
$I_{(m)}=I'_{(m)}$. This defines an equivalence relation $\sim$ 
in $\sigmn$. Lemma \ref{lemm:modmface} implies that,
if $X$ is an equivalence class, the members of $X$ have some $K$ as 
the common $\bmod\: m$ face. 
We shall call this $K$ the \emph{core} of the equivalence class $X$.

\begin{lemm}\label{lemm:type}
Suppose that the condition (P) is satisfied. Let $X$ be an equivalence 
class of $(v,m)$-equivalence relation.
For $x\in H_T^2(\Delta)$ the value 
$\langle \iota_I^*(x),v\rangle \bmod\: m$
does not depend on the choice of $I$ in $X$. 
\end{lemm}
\begin{proof}
Write $x=\sum_{i\in \sigmone}a_ix_i$. Let $K$ be the core of $X$. Then
\[ \langle \iota_I^*(x),v\rangle =
 \sum_{i\in K}a_i\langle u_i^I,v\rangle
  +\sum_{i\in I,\notin K}a_i\langle u_i^I,v\rangle. \]
Since $\langle u_i^I,v\rangle =0\mod\: m$ for $i\notin K$ 
by the definition of the core and
$\langle u_i^I,v\rangle \mod\: m$ for $i\in K$ does not depend on 
$I$ in $X$ by Lemma \ref{lemm:modmface}, $\langle \iota_I^*(x),v\rangle \mod m$
does not depend on the choice of $I$ in $X$.   
\end{proof}

\begin{coro}\label{coro:type}
Suppose that the condition (P) is satisfied. 
Assume that $c_1(\Delta,\V)$ is divisible by $N$, and write
$c_1^T(\Delta,\V)=Nx+u, u\in H^2(BT)$. Let $v$ and $m$ be as above.
If we write $\langle u_i^I,v\rangle$ in the form
\[ \langle u_i^I,v\rangle =mh_i +r_i \ \text{with}\ 0\leq r_i<m ,\]
then the sum $\sum_{i\in I}h_i\bmod\: N$ depends only on the 
$(v,m)$-equivalence class $X$ of $I$.
\end{coro}
\begin{proof}
Let $K$ be the core of $X$. 
Then $\l u_i^I, v\r=0 \bmod m$ for $i\in I\setminus K$ and 
$\l u_i^I, v\r \bmod m$ for $i\in K$ is independent of $I$ in $X$ 
by Lemma \ref{lemm:modmface} so that $r_i=0$ for $i\in I\setminus K$ 
and $r_i$ for $i\in K$ is independent of $I$. Hence 
the sum $\sum_{i\in I}r_i$ is a constant depending only on $X$ 
which we shall denote by $r$.
We put $h_I=\sum_{i\in I}h_i$. Then
$\langle u^I,v\rangle=mh_I+r$.
By Lemma \ref{lemm:type} $\langle \iota_I^*(x),v\rangle$ is of the form
\[ \langle \iota_I^*(x),v\rangle =mh'_I +r' ,\]
for $I\in X$, where $r'$ is independent of $I$. 
Therefore, if we write $\langle u,v\rangle =r''$, then
\[ \langle u^I,v\rangle =Nmh'_I +Nr' +r'' .\]
If we compare this with
\[ \langle u^I,v\rangle =mh_I +r ,\]
we see that $Nr'+r''$ is of the form $Nr'+r''=mh'+r$ and
$h_I=Nh'_I+h'$. 
This shows that $h_I\bmod\,N$ depends only on $X$.
\end{proof}
Under the situation of Corollary \ref{coro:type}, the $\bmod\,N$
value of $\sum_{i\in I}h_i$ will be called $(v,m)$-\emph{type} of 
$X$ and will be 
denoted by $h(v,m,X)$. Similarly the $\bmod\,N$
value of $\langle u^I,v\rangle$ (which is independent of 
$I\in \sigmn$ by Lemma \ref{lemm:Ndivisible})
will be called $v$-\emph{type} and will be 
denoted by $h(v)$.

\begin{lemm}\label{lemm:many types}
Suppose that the condition (P) is satisfied. 
Assume $c_1(\Delta,\V)$ is divisible by $N$. Any non-zero 
$b\in\Z/N$ can occur as $v$-type when $v$ varies over generic
vectors in $L_\V$.
\end{lemm}
\begin{proof}
This follows readily from the fact that $\{v_i|i\in I\}$ 
is a basis of $L_\V$ for each $I\in \sigmn$.
\end{proof}


\section{Rigidity theorem}\label{sec:rigidity}

Let $\Delta=(\Sigma,C,w^\pm)$ be a complete simplicial multi-fan and 
$\V$ a set of prescribed edge vectors. Let $N$ be an integer greater 
than $1$. When $\zeta^N=1,\zeta\not=1$,  $\var(\Delta,\V)$ and 
$\hatvar(\Delta,\V)$ are called elliptic genus of \emph{level} $N$
and {orbifold elliptic genus of \emph{level} $N$ respectively. 
Let $v\in L_\V$ be a generic vector. 
Since $\zeta^N=1$, $\varv(\Delta,\V)$ and $\hatvarv(\Delta,\V)$ are 
elliptic functions in $z$ with respect to the lattice 
$\Z N\tau \oplus \Z$, 
because $\phi(z,\tau,\sigma)$ is such a function by \eqref{eq:transformula2}.

Hereafter the condition (P) is assumed 
throughout this section.
\begin{lemm}\label{lemm:modNtransform}
Suppose that the condition (P) is satisfied and $c_1(\Delta,\V)$ is 
divisible by $N$. Let $v\in L_\V$
be a generic vector and $h(v)$ the $v$-type. Then the elliptic
genus $\varv(\Delta,\V)$ along $v$ and orbifold elliptic genus 
$\hatvarv(\Delta,\V)$ along $v$  of level $N$ transform by
\begin{align*}
 \varv (\Delta,\V;z+\tau, \tau,\sigma)&
=\zeta^{h(v)}\varv(\Delta,\V;z,\tau,\sigma) ,\\
 \hatvarv(\Delta,\V;z+\tau,\tau,\sigma)&
=\zeta^{h(v)}\hatvarv(\Delta,\V;z,\tau,\sigma).
\end{align*} 
\end{lemm}
\begin{proof}
\[ \begin{split}
 \prod_{i\in I}\phi(\l u_i^I,-(z+\tau)v+&\tau v(h_1)-v(h_2)\r,\tau,\sigma) \\
 &= \zeta^{\sum_{i\in I}\langle u_i^I,v\rangle}
 \prod_{i\in I}\phi(\l u_i^I,-zv+\tau v(h_1)-v(h_2)\r,\tau,\sigma)
 \end{split} \]
by \eqref{eq:transformula2}. But $\sum_{i\in I}\langle u_i^I,v\rangle
 =h(v) \mod N$ which is independent of $I$. Hence we obtain
\[
 \varv (\Delta,\V;z+\tau, \tau,\sigma)
=\zeta^{h(v)}\varv(\Delta,\V;z,\tau,\sigma),\ \ 
 \hatvarv(\Delta,\V;z+\tau,\tau,\sigma)
=\zeta^{h(v)}\hatvarv(\Delta,\V;z,\tau,\sigma),
\] 
since $\zeta^N=1$. 
\end{proof}

The following theorem and corollary are versions of rigidity 
theorem and vanishing theorem for multi-fans.
\begin{theo}\label{theo:rigid}
Let $\Delta$ be a complete simplicial multi-fan and $\V$ a set of 
prescribed edge vectors satisfying the condition (P). Let
$v\in L_\V$ be a generic vector. Assume that $c_1(\Delta,\V)$ 
is divisible by an integer $N>1$. Then the equivariant elliptic genus 
$\varv(\Delta,\V;z,\tau,\sigma)$ of level $N$ along $v$
is rigid, i.e. it is constant as a function of $z$.
\end{theo}
\begin{coro}\label{coro:vanishing}
Under the same situation as in Theorem \ref{theo:rigid} the 
equivariant elliptic genus $\var(\Delta,\V)$ of level $N$
constantly vanishes. In particular, if 
$\Delta$ is a complete non-singular multi-fan whose first Chern class 
$c_1(\Delta)$ is divible by $N$, then its equivariant elliptic 
genus $\var(\Delta)$ of level $N$ constantly vanishes.
\end{coro}
\begin{proof}
We postpone the proof of Theorem \ref{theo:rigid}. As to 
Corollary \ref{coro:vanishing}, we take a generic vector
$v$ such that $\zeta^{h(v)}\not=1$, which is possible by
Lemma \ref{lemm:many types}. Since $\varv(\Delta,\V;z,\tau,\sigma)$ 
is constant by
Theorem \ref{theo:rigid} for any $v$, $\var(\Delta,\V)$ is also
constant, which we denote by $\var$. Since $\var =
\zeta^{h(v)}\var$ by 
Lemma \ref{lemm:modNtransform}, 
$\var$ must be equal to $0$. 
\end{proof}

The degree $0$ term in the $q$-expansion of $\var(\Delta,\V)$ 
reduces to the $T_y$-genus $T_y[\Delta]$. It is independent 
of $\V$ as was remarked in Section 3. We obtain
\begin{coro}\label{coro:Tyvanish}
Suppose that $\Delta$ is a complete simplicial multi-fan and 
there is a set of generating vectors $\V$ which satisfies 
the condition (P). If $c_1(\Delta,\V)$ is 
divisible by $N$, then the $T_y$-genus
$T_y[\Delta]$ of $\Delta$ vanishes for $(-y)^N=1,\ -y\not=1$.
\end{coro}
 
If $M$ is a stably almost complex closed manifold, 
then the signature $\Sign(M)$ equals $T_1(M)$. After this 
the $T_y$-genus for $y=1$ of a complete non-singular multi-fan
$\Delta$ is called the signature of $\Delta$ and is denoted by 
$\Sign(\Delta)$. 
Hence
\begin{coro}\label{coro:spin}
The signature $\Sign(\Delta)$ of a complete non-singular
multi-fan $\Delta$ with $c_1(\Delta)=0 \bmod\: 2$ vanishes.
\end{coro}

\begin{rema}
If $M$ is a torus manifold, its elliptic genus coincides with
that of its associated multi-fan $\Delta(M)$. If $c_1(M)$ is 
divisible by $N$, $c_1(\Delta(M))$ is also divisible by $N$ 
as was remarked in Section \ref{sec:chern}, and hence its equivariant 
elliptic genus of level $N$ vanishes. In case
$N=2$ we have the following conclusion. 
The equivariant Stiefel-Whitney class $w_2^T(M)\in H_T^2(M;\Z/2)$ is 
defined and descends to 
$c_1^T(\Delta(M)) \bmod\,2$. If $M$ is a spin torus manifold, 
then $w_2^T(M)$ lies in $H^2(BT;\Z/2)$. Therefore 
$c_1^T(\Delta(M))$ is divisible by $2$. It follows from
Corollary \ref{coro:spin} that the signature $\Sign(M)$ of 
$M$ vanishes. This can be also deduced from Corollary in
1.5 of \cite{HS}. A complete non-singular multi-fan with
$c_1(\Delta)=0 \bmod\: 2$ might be called a spin multi-fan.
\end{rema}

The rest of this section is devoted to the proof of Theorem
\ref{theo:rigid}. It is convenient to consider 
$\varv (\Delta,\V;z, \tau,\sigma)$ as a function of 
$t=e^{2\pi\img z}$ rather than $z$. It will be denoted by 
$\varv(t)$. 

We assume throughout that $c_1(\Delta,\V)$ is 
divisible by $N$, and $\sigma=\frac{k}{N},\ 0<k<N$, so that 
$\zeta^N=1$. It suffices to show that 
$\varv(t)$ has no poles since it is an elliptic function.
It is clear that possible poles $t=\lam$ satisfy $\lam^m q^{-s}=1$
for some integers $m\geq 1$ and $s$. Hence it suffices to
show that $\varv(tq^{\frac{s}{m}})$ has no poles $t=\lam$ with
$\lam^m=1$. Let 
\[ \sigmn=\sqcup_\nu \Sigma_\nu^{(n)} \]
be the decomposition into $(v,m)$-equivalence classes. We fix a class 
$\Sigma_\nu^{(n)}$ 
and write
\[ \langle u_i^I,v\rangle =m_i =mh_i+r_i,\ 0\leq r_i<m, \]
for $I\in \Sigma_\nu^{(n)}$. Let
$K_\nu$ denote the core of $\Sigma_\nu^{(n)}$.
Remember that 
\[ r_i=0 \ \Leftrightarrow \ i\notin K_\nu ,\]
and $r_i$ is independent of $I\in\Sigma_\nu^{(n)}$ for $i\in K_\nu$. 
Note that $\sum_{i\in I}h_i =h(v,m,\Sigma_\nu^{(n)}) \mod N$.
We set
\[ \bar{v}=\sum_{i\in K}r_iv_i. \]
Then $\l u_i^I,v-\bar{v}\r=mh_i$. 
A straightforward calculation using \eqref{eq:transformula2} yields
\[ \varv(tq^{\frac{s}{m}})=\sum_\nu\zeta^{sh(v,m,\Sigma_\nu^{(n)})}
     \var_\nu(t) ,\]
where            
\[
 \begin{split}
  \var_\nu(t) =& 
 \sum_{I\in \Sigma_\nu^{(n)}}\frac{w(I)}{|H_I|}\cdot \\
   &\sum_{h\in H_I}
    \prod_{i\in I\setminus K_\nu}
    \phi(-\l u_i^I, zv+v(h)\r, \tau, \sigma)
    \prod_{i\in K_\nu}
    \phi(-\l u_i^I, zv+\frac{s}{m}\tau\bar{v}+v(h)\r, \tau, \sigma).
 \end{split}
\]

\begin{lemm}\label{lemm:vmlaurent}
Fix $\nu$. Then $\var_\nu(t)$ can be expanded in the following form:
 \[ \var_\nu(t)=
    \sum_{n=0}^\infty \bar{\var}_n(t)q^{\frac{n}{m}} \]
with $\zeta^{\frac{n}{2}}\bar{\var}_n(t)\in R(S^1)\otimes 
\Z[\zeta,\zeta^{-1}]$.
\end{lemm}
\begin{proof}
The proof is similar to that of Proposition \ref{prop:laurent}.
Put $\tau_1=\frac{s}{m}\tau$. We expand $\var_\nu(t)$ in the form
\[ \var_\nu(t)= 
 \sum_{s_1\in \Z,s_2\in \Z_{\geq 0}}\var_{s_1,s_2}(z,q_1)\zeta^{s_1}q^{s_2}. \]
Then 
$\zeta^{\frac{n}{2}}\var_{s_1,s_2}(z,q_1)\zeta^{s_1}q^{s_2}$ is a 
sum of expressions of 
the following form.
\begin{equation*}\label{eq:qone2}
 \sum_{I\supset K_\nu}\frac{w(I)}{|H_I|}
 \sum_{l\in \alpha_I} 
  \sum_{h\in H_I}
 \frac{\prod_{i\in I}\chi_I(l,h)t^{\l l,v\r}
 q_1^{-\l l,\bar{v}\r}}
  {\prod_{i\in I\setminus K_\nu}(1-\chi_I(u_i^I,h)^{-1}t^{-\l u_i^I,v\r})}, 
\end{equation*}
where $\alpha_I\in \mathcal{O}_I$ satisfying 
$\tau_{I,I'}(\alpha_I)=\alpha_{I'}$.
It follows that $\zeta^{\frac{n}{2}}\var_{s_1,s_2}(z,q_1)\zeta^{s_1}q^{s_2}$ 
belongs to $\Z[t,t^{-1}, q_1,q_1^{-1}]$ by 
Lemma \ref{lemm:integrality}. Noting that $q_1=q^{\frac{s}{m}}$, we see 
that $\var_\nu(t)$ can be expanded in the form
 \[ \var_\nu(t)=
    \sum_{n=0}^\infty \bar{\var}_n(t)q^{\frac{n}{m}} \]
with $\zeta^{\frac{n}{2}}\bar{\var}_n(t)\in R(S^1)\otimes 
\Z[\zeta,\zeta^{-1}]$.
\end{proof}

We shall show that each $\var_\nu(t)$ has no 
poles at $\lambda$ with $\lam^m=1$. Suppose that
$\lam$ is a possible pole with $\lam^m=1$.
Let $\bar{\var}_{n,I}(t)$ be the contribution from $I$ in 
$\bar{\var}_n(t)$. There is an open set $U$ containing $\lam$
such that $\bar{\var}_{n,I}(t)$ is holomorphic in 
$U\setminus \{\lam\}$ for each $I$. 
$\sum_{n=0}^\infty\bar{\var}_{n,I}(t)q^{\frac{n}{m}}$ converges uniformly 
in any compact set in $U-\{\lam\}$. Note that 
$\sum_I\bar{\var}_{n,I}(t)=\bar{\var}_n(t)$ is holomorphic in
$U$ because it is a finite Laurent series by Lemma 
\ref{lemm:vmlaurent}. 
We now quote a lemma from \cite{H}.
\begin{lemm}[Hirzebruch]\label{lemm:holo}
Let $\{b_{n,j}\}$ be a family of meromorphic functions on $U$ with $n$ running 
over non-negative integers and $j$ running over 
some finite set $J$. Suppose that they satisfy the following 
properties:
\begin{enumerate}
 \item $b_{n,j}$ is holomorphic in $U\setminus \{\lam\}$,
 \item $b_n=\sum_jb_{n,j}$ is holomorphic in $U$,
 \item $\sum_{n=0}^\infty b_{n,j}$ converges uniformly in any compact
set in $U\setminus \{\lam\}$.
\end{enumerate}
Then $\sum_{n=0}^\infty b_n$ converges uniformly in any compact
set in $U$ and is a holomophic extension of 
$\sum_{j\in J}\sum_{n=0}^\infty b_{n,j}\mid U\setminus\{\lam\}$.
\end{lemm}
We apply this Lemma to $\{\bar{\var}_{n,I}(t)q^{\frac{n}{m}}\}$. 
It follows that $\var_\nu(t)$
and hence $\varv(tq^{\frac{s}{m}})$ has no pole
at $t=\lam$. Hence $\varv(t)$ has no poles and it is a
constant. This proves Theorem \ref{theo:rigid}.

\begin{rema} 
Suppose that the pair $(\Delta,\V)$ satisfies the condition (P).
In Remark \ref{rem:cover} we introduced a non-singular multi-fan 
$\tilde{\Delta}$ in a lattice $\tilde{L}$ and an injective homomorphism
$\pi:\tilde{L}\to L$ with image $\pi(\tilde{L})=L_\V$. 
We identify $L^*$ with the sublattice $\pi^*(L^*)$ of $\tilde{L^*}$.
If we write $\var(\tilde{\Delta})$ as
\[ \var(\tilde{\Delta})=\sum_{u\in \tilde{L}^*}a_ut^u \]
the coefficients $a_u$ are given by Theorem \ref{thm:BL}.
Comparing this with the expression of $\var(\Delta,\V)$ 
given by Theorem \ref{thm:BL} and noting
that $\langle u,\tilde{v}_i\rangle =\langle u,v_i\rangle$
we see that 
\begin{equation*}\label{eq:varyn}
 \var(\Delta,\V)=\sum_{u\in L^*}a_ut^u .
\end{equation*} 
Thus, if the coefficients of $t^u$ in $\var(\tilde{\Delta})$ 
vanish for all $u\in \tilde{L}^*$, then those of $\var(\Delta,\V)$ 
also vanish. This is the case when $c_1(\Delta,\V)$ is divisible by $N$ 
and $\zeta^N=1$, since
$c_1(\tilde{\Delta})$ is also divisible by $N$, 
as was remarked in Remark \ref{rem:N-divisible}.
\end{rema}

\begin{rema}
Examples show that the rigidity phenomenon does not necessarily 
occur for pairs $(\Delta,\V)$ without condition (P). 
\end{rema}


\section{Applications}
Let $\Delta$ be a complete simplicial multi-fan of dimension $n$.
\begin{lemm}
If the Todd genus $T_0[\Delta]$ equals $1$ and $w(I)=1$ for all 
$I\in \sigmn$, then 
\begin{equation}\label{eq:hq}
 h_k(\Delta)=\#\{I\in \sigmn \mid \mu(I)=k\},
\end{equation}
in \eqref{eq:h} and 
\begin{equation}\label{eq:ek}
  e_k(\Delta)=\#\sigmk 
\end{equation}
in (\ref{eq:T_y}).
\end{lemm}
\begin{proof}
\eqref{eq:hq} is immediate. As was remarked in Note after 
Proposition \ref{prop:Tygenus}, $T_0(\Delta)$ equals $\deg(\Delta)$. 
Therefore, if $T_0(\Delta)=1$ then $\deg(\Delta)=1$, and consequently 
$\deg(\Delta_J)=1$ for all $J\in \Sigma$, as is easily seen. 
Then \eqref{eq:ek} follows.
\end{proof}
\begin{note}
Note that these conditions are always satisfied for complete simplicial
ordinary fans. In particular, if $M$ is a complete non-singular toric variety, 
then the fan $\Delta(M)$ associated with $M$ satisfies these conditions.
It is known that $T_{-t^2}[\Delta(M)]$ is equal to the Poincar\'e 
polynomial $P(t)$ of $M$ if $M$ is a non-singular projective toric 
variety, see e.g. \cite{F}.
\end{note}

\begin{prop}\label{prop:N=n+1}
Let $\Delta$ be a complete non-singular multi-fan with Todd
genus $T_0[\Delta]\not=0$. If 
$c_1(\Delta)$ is divisible by a positive integer $N$, then 
$N$ is equal to or less than $n+1$. In the cases $N=n+1$ and $N=n$
the $T_y$- genus must be of the following forms
\begin{equation}\label{eq:Tyn+1}
 T_y[\Delta]=T_0[\Delta]\sum_{k=0}^n(-y)^k \quad (N=n+1),
\end{equation} 
and
\begin{equation}\label{eq:Tyn}
 T_y[\Delta]=T_0[\Delta](1-y)\sum_{k=0}^{n-1}(-y)^k \quad (N=n).
\end {equation}
\end{prop}
\begin{proof}
Suppose that $c_1(\Delta)$ is divisible by $N$. Then, 
by Corollary \ref{coro:Tyvanish} $T_y[\Delta]$ considered as
a polynomial in $-y$ has roots at all $N$-th roots of unity
other than $1$. Hence it must be divisible by $\sum_{k=0}^{N-1}(-y)^k$.
On the other hand it is a polynomial of degree $n$ with constant
term $T_0[\Delta]\not=0$ by Proposition
\ref{prop:Tygenus}. Therefore we must have $N-1\leq n$.

Suppose that $N=n+1$. Then the same reasoning as above proves
(\ref{eq:Tyn+1}).
If $N=n$, then $T_y[\Delta]$ is divisible by 
$\sum_{k=0}^{n-1}(-y)^k $. Since the constant term and the 
coefficient (as a polynomial of $-y$) of the highest
term do not vanish by assumption and Proposition \ref{prop:Tygenus},
$T_y[\Delta]$ must be of the form (\ref{eq:Tyn}).
\end{proof}
\begin{lemm}\label{lemm:N=n+1 and n}
Let $\Delta$ be a complete non-singular multi-fan with 
$T_0[\Delta]=1$ and $w(I)=1$ for all $I\in \sigmn$. If $T_y[\Delta]$
is of the form (\ref{eq:Tyn+1}), then
\[ \#\sigmone =n+1 \ \text{and}\ \#\sigmn =n+1 .\]
If $T_y[\Delta]$ is of the form (\ref{eq:Tyn}), then
\[ \#\sigmone =n+2 \ \text{and}\ \#\sigmn =2n .\] 
Moreover, in case $n\geq 3$, $\#\Sigma^{(2)} =\frac12n(n+3)$.
\end{lemm}
\begin{proof}
The equality (\ref{eq:Tyn+1}) with $T_0[\Delta]=1$ implies that 
$h_k(\Delta)=1$ for all $k$ with $0\leq k\leq n$ by 
Proposition \ref{prop:Tygenus}. 
This implies that $\#\sigmn =n+1$ by (\ref{eq:hq}). 
Then, using (\ref{eq:T_y}) and (\ref{eq:ek}) we see 
that $\#\sigmone =e_1(\Delta)=n+1$.

Similarly the equality (\ref{eq:Tyn}) with $T_0[\Delta]=1$
implies that $h_k(\Delta)=1$ for $k=0,n$ and $h_k(\Delta)=2$ 
for $1\leq k\leq n-1$. 
This implies that $\#\sigmn =2n$ by (\ref{eq:hq}), 
and yields, together with 
(\ref{eq:T_y}) and (\ref{eq:ek}), 
the equalities 
$\#\sigmone =e_1(\Delta)=n+2$ and $\#\Sigma^{(2)} =e_2(\Delta)=\frac12n(n+3)$
in case $n\geq 3$.
\end{proof} 

\begin{coro}\label{coro:projective space}
Let $M$ be a complete non-singular toric variety of dimension $n$. If 
$c_1(M)$ is divisible by $n+1$, then $M$ is isomorphic to 
the projective space $\Proj^n$ as a toric variety.
\end{coro}
\begin{proof}
By Proposition \ref{prop:N=n+1} and Lemma \ref{lemm:N=n+1 and n}
the fan associated with $M$ has $n+1$ $1$-dimensional cones and
$n+1$ $n$-dimensional cones. Such a fan is unique (up to automorphisms
of the lattice $L$) and coincides with the fan associated with 
$\Proj^n$. Since a toric variety is determined by its fan,
$M$ must be $\Proj^n$.
\end{proof}

In order to handle the case $N=n$ we investigate the fan associated
with a projective space bundle over a projective space. Let
$\xi$ denote the hyperplane bundle (dual of the 
tautological line bundle) over $\Proj^r$. Let $1\leq r<n$ and set
$\eta =(\bigoplus_{i=r+1}^n\xi^{k_i})\oplus 1$ where $k_i$ are integers
and $1$ denotes the trivial line bundle. The associated 
projective space bundle of $\eta$ will be denoted by $M$. 
It is a complex manifold. A point of $M$ is expressed
in homogeneous coordinate
\begin{equation}\label{eq:projective eta}
  [z_0,z_1,\ldots ,z_r,w_{r+1},\ldots ,w_n, w_{n+1}] 
\end{equation}
where $z_i, w_j \in \C,\ (z_0,z_1,\ldots ,z_r)\not= (0,0,\ldots,0),
\ (w_{r+1},\ldots ,w_n,w_{n+1})\not= (0,\ldots ,0,0) $, and
if, $\alpha\in \C^*$, then 
\[ [\alpha z_0,\alpha z_1,\ldots ,\alpha z_r,
 \alpha^{k_{r+1}}w_{r+1},\ldots ,\alpha^{k_n}w_n, w_{n+1}] \]
and
\[ [z_0,z_1,\ldots ,z_r,
     \alpha w_{r+1},\ldots ,\alpha w_n,\alpha w_{n+1}] \]
are identified with (\ref{eq:projective eta}). 

Let the $(n+1)$-dimensional torus $T^{n+1}=S^1\times \dots \times S^1$ 
act on $M$ by
\begin{equation*}\label{eq:Taction}
 \begin{split}
(t_0,t_1,\ldots ,t_n)
      &[z_0,z_1,\ldots ,z_r,w_{r+1},\ldots ,w_n, w_{n+1}] \\
&=[t_0z_0,t_1z_1,\ldots ,t_rz_r,t_{r+1}w_{r+1},\ldots ,t_nw_n, w_{n+1}]
 \end{split}
\end{equation*}
The action is a holomorphic action. The subgroup 
$D'=\{(t,\ldots ,t,t^{k_{r+1}},\ldots ,t^{k_n})\mid t\in S^1\}$
of $T^{n+1}$ acts trivially on $M$. Hence the quotient $T=T^{n+1}/D'$
acts on $M$. Put
\begin{equation*}
 M_i= \left\{
  \begin{array}{ll}
     \{[z_0,z_1,\ldots ,z_r,w_{r+1},\ldots ,w_n, w_{n+1}]
       \mid z_i=0 \} & \text{for} \ \ 0\leq i\leq r  \\
     \{[z_0,z_1,\ldots ,z_r,w_{r+1},\ldots ,w_n, w_{n+1}]
       \mid w_i=0 \} & \text{for} \ \ r+1\leq i\leq n+1 .
  \end{array}
 \right.
\end{equation*} 
We also put
\begin{align*}
 S_i &= \{(1,\ldots ,1,t_i,1,\ldots ,1)\in T^{n+1} \}
    \ \  \text{for} \ \ 0\leq i\leq n  \\
 \intertext{and}
 S_{n+1} &=\{(1,1,\ldots ,1,t_{r+1},\ldots ,t_n)\in T^{n+1}
         \mid t_{r+1}= \dots =t_n \} .
\end{align*}
We shall denote by the same letter the image of $S_i$ 
in $T$. It is easy to see that $S_i$ pointwise fixes 
$M_i$, and there are no other circle subgroups of $T$ which
have (complex) codimension $1$ fixed point set components. 

Let $\tilde{v}_i\in \Hom(S^1,T^{n+1})$ denote the inclusion 
homomorphism of $S^1$ into the $i$-th factor of $T^{n+1}$ for 
$0\leq i\leq n$. Put 
$v_i=\pi_*(\tilde{v}_i)\in \Hom(S^1,T)$ where $\pi_*$ is the
homomophism induced by the projection 
$\pi:T^{n+1}\to T$. From the definition it follows that there is a relation
\begin{equation}\label{eq:vzero}
 v_0+v_1+\cdots +v_r +k_{r+1}v_{r+1}+\cdots +k_nv_n =0 .
\end{equation} 
We also put 
\begin{equation}\label{eq:vn+1}
v_{n+1}=-(v_{r+1}+\dots +v_n)\in \Hom(S^1,T). 
\end{equation}
Then $S_i$ is the image of $v_i:S^1\to T$. Moreover $S^1$ acts 
via $v_i$ on
each fiber of the normal bundle of $M_i$ in $M$ by standard
complex multiplication. Thus $\{v_i\mid i=0,\ldots ,n,n+1\}$ 
coincides with the set of primitive edge vectors of 
$1$-dimensional cones of the multi-fan $\DeltaM=(\SigmaM,C(M),w(M)^\pm)$ 
associated with the torus manifold $M$, and we have 
$\SigmaM^{(1)}=\{0,1,\ldots,n,n+1\}$, cf. \cite{HM}.

To determine the whole augmented simplicial set $\SigmaM$, we need to
look at the fixed point set $M^T$. For $i\in \{0,1,\dots r\}$, put
$I_i=\{0,1,\dots r\}\setminus \{i\}$, and for $j\in \{r+1,\dots n+1\}$,
put $J_j=\{r+1,\dots n+1\}\setminus \{j\}$.
It is not difficult to see that $M^T$ consists of points
\[ M_{I_i}\cap M_{J_j} , \ \ i\in \{0,1,\dots r\},
            \ \ j\in \{r+1,\dots n+1\},\]
where $M_{I_i}=\cap_{k\in I_i}M_k$ and $M_{J_j}=\cap_{l\in J_j}M_l$.
This implies that 
\begin{equation*}\label{eq:sigmn}
 \SigmaM^{(n)}=\{I_i\cup J_j\mid 
     i\in \{0,1,\dots r\},\ j\in \{r+1,\dots n+1\} \}.
\end{equation*}
In particular 
\begin{equation*}\label{eq:number of sigmn}
 \#\SigmaM^{(n)}=(r+1)(n-r+1) .
\end{equation*} 
It follows that
\begin{equation*}\label{eq:minimum} 
 \#\SigmaM^{(n)}\geq 2n \ \ \text{and}\ \ 
 \#\SigmaM^{(n)}=2n \ \ \text{if and only if}\ \ 
            r=1 \ \ \text{or}\ \  n-1.
\end{equation*} 
     
Let $\tau$ be the tautological line bundle over the projective
space bundle $M$. Its dual $\tau^*$ restricts to the hyperplane bundle 
on each fiber of $\pi:M\to \Proj^r$.
Let $\omega\in H^2(M)$ be the first Chern class of $\tau^*$.
Then, by the Leray-Hirsch theorem, $H^*(M)$ is a free 
$H^*(\Proj^r)$-module on generators 
$1,\omega, \omega^2,\ldots ,\omega^{n-r}$. In particular,
$H^2(M)$ is a free module on $\omega, \omega'$, where $\omega'$
is the image of the canonical generator of $H^2(\Proj^r)$
by $\pi^*$. We have
\begin{lemm}\label{lemm:chern class}
\[ c_1(M)=(n-r+1)\omega +(\sum_{i=r+1}^nk_i +r+1)\omega' .\]
\end{lemm}
\begin{proof}
The tautological line bundle $\tau$ is a subbundle of $\pi^*\eta$,
and the tangent bundle along the fibers $T_fM$ of 
$\pi:M\to \Proj^r$ is isomorphic to $\Hom(\tau,\pi^*\eta/\tau)=
\tau^*\otimes (\pi^*\eta/\tau)$. Hence
\begin{align*}
 c_1(T_fM) &=(n-r)c_1(\tau^*)+c_1(\pi^*\eta/\tau) \\
             &=(n-r)\omega+c_1(\pi^*\eta)-c_1(\tau) \\
             &=(n-r)\omega+(\sum_ik_i)\omega'+\omega \\
             &=(n-r+1)\omega+(\sum_ik_i)\omega'.
\end{align*}
Since the tangent bundle $TM$ is isomorphic to 
$\pi^*T\Proj^r\oplus T_fM$, and $c_1(\pi^*T\Proj^r)=(r+1)\omega'$,
we have
\[ c_1(M)=(n-r+1)\omega+(\sum_ik_i+r+1)\omega'. \]
\end{proof}
As an immediate consequence of Lemma \ref{lemm:chern class} we obtain
\begin{coro}\label{coro:bundle over Pone}
Let $M=\Proj(\eta)$ be as above. Then $c_1(M)$ is divisible by $n$
if and only if $r=1$ and $\sum_{i=r+1}^nk_i+2$ is divisible by $n$.
\end{coro}

We now consider complete non-singular multi-fans having first Chern 
class divisible by $n$.
\begin{lemm}\label{lemm:modn}
Let $\Delta=(\Sigma,C,w^\pm)$ be a complete non-singular multi-fan 
of dimension $n$ such that
\[ T_0[\Delta]=1,\  w(I)=1 \ \ \text{for all}\ \ I\in \sigmn, \]
\[  \#\sigmone=n+2,\ \  \#\sigmn=2n,\ \ \text{and}\ \ 
  \#\Sigma^{(2)}=\frac12n(n+3)\ \ \text{in case}\ \ n\geq 3. \]
Then it is equivalent to the fan associated to a $\Proj^{n-1}$  
bundle over $\Proj^1$ or a $\Proj^1$-bundle over 
$\Proj^{n-1}$.  
\end{lemm}
\begin{proof}
Let $\{v_i\}_{i=0}^{n+1}$ be the primitive edge vectors of the 
$1$-dimensional cones. In view of (\ref{eq:vzero}) and 
(\ref{eq:vn+1}) it suffices to show that, under a suitable 
numbering, they satisfy the relations
\begin{equation}\label{eq:vPone}
 v_1+v_2+\cdots +v_n=0,\ \ v_0+v_{n+1}+\sum_{i=2}^nk_iv_i=0 ,
\end{equation} 
or
\begin{equation}\label{eq:vn-1}
 v_1+v_2+\cdots +v_n+kv_{n+1}=0,\ \ v_0+v_{n+1}=0 .
\end{equation}
We first deal with the case $n\geq 3$. From the completeness
we see that each $1$-dimensional cone is a face of at least
$n$ $2$-dimensional cones, and it is a face of at most $n+1$
$2$-dimensional cones because the number of $1$-dimensional
cones are $n+2$. Since the number of $2$-dimensional cones
is $\frac12n(n+3)=n+\frac12n(n+1)$, we conclude that there are two
edge vectors, say $v_0$ and $v_{n+1}$ such that $v_0$ 
$(v_{n+1}\ \ \text{respectively})$ 
spans $2$-dimensional cones with each of remaining vectors
$v_1,v_2,\ldots ,v_n$, and each $v_i,\ 1\leq i\leq n,$
spans $2$-dimensional cones with $v_j,\ j\not=i$. Thus the projected
multi-fan $\Delta_{\{0\}}$ has exactly $n$ $1$-dimensional cones.
It is complete and non-singular as a projected multi-fan of a complete 
non-singular multi-fan $\Delta$. It follows that $\Delta_{\{0\}}$ 
is equivalent to
the fan of $\Proj^{n-1}$, and the projected edge vectors
$\bar{v}_i,\ 1\leq i\leq n,$ satisfy the relation
\[ \bar{v}_1+\bar{v}_2+\cdots +\bar{v}_n=0 .\] 
This implies the relation 
\[ v_1+v_2+\cdots +v_n=kv_0 .\]
Similarly we have      
\[ v_1+v_2+\cdots +v_n=k'v_{n+1} .\]
If $k=0$ then $k'=0$ since $v_{n+1}\not=0$, and $v_1,v_2,\ldots ,v_n$
lie on a hyperplane $v_1+v_2+\cdots +v_n=0$. Since 
the multi-fan $\Delta$ is complete and non-singular, the primitive
vectors $v_0$ and $v_{n+1}$ lie on the different sides of that 
hyperplane and must satisfy a relation as described in (\ref{eq:vPone}).

If $k\not=0$ then $k'\not=0$ and $v_0$ and $v_{n+1}$ are linearly
dependent primitive vectors. Therefore we must have $v_0+v_{n+1}=0$.
Thus (\ref{eq:vn-1}) holds. This proves Lemma \ref{lemm:modn}
in the case $n\geq 3$.

The case $n=2$ is similar and easier. We see that there are four 
primitive edge vectors $v_0,v_1,v_2,v_3$ in $2$-dimensional
vector space $V=L\otimes \R$ such that
\[ v_1 +v_2=kv_0=k'v_3. \]
By the same reasoning as above we derive 
\[  v_1 +v_2=0,\ v_0+v_3 +kv_2=0 \ \ \text{or}\ \ 
    v_1 +v_2+kv_3=0, \ v_0+v_3=0 .\]
\end{proof}
\begin{coro}\label{coro:projective bundle}
Let $M$ be a complete non-singular toric variety of dimension $n$. If 
$c_1(M)$ is divisible by $n$, then $M$ is isomorphic to 
an $(n-1)$-dimensional projective space bundle $\Proj(\eta)$ 
with $\eta=(\bigoplus_{i=2}^{n}\xi^{k_i})\oplus 1$ over $\Proj^1$ 
such that $(\sum_{i=2}^{n}k_i)+2$ is divisible by $n$. One dimensional 
cones of the fan associated with $M$ is generated by primitive 
vectors $v_0,v_1,v_2,\ldots,v_n,v_{n+1}$ satisfying 
\[ v_0+v_{n+1}+\sum_{i=2}^nk_iv_i=0,\ v_1+v_2+\cdots +v_n=0. \]
\end{coro}
\begin{proof}
By Proposition \ref{prop:N=n+1} and Lemma \ref{lemm:N=n+1 and n}
the fan $\Delta(M)$ associated with $M$ has $n+2$ $1$-dimensionl 
cones and $2n$ $n$-dimensional cones. Moreover the number of 
$2$-dimensional cones is $\frac12n(n+3)$ in case $n\geq 3$. By 
Lemma \ref{lemm:modn} $\Delta(M)$ is equivalent to that of a
$\Proj^{n-1}$-bundle over $\Proj^1$ or $\Proj^1$-
bundle over $\Proj^{n-1}$. Among such manifolds those with 
$c_1$ divisible by $n$ are of the form given in Corollary 
\ref{coro:bundle over Pone}.
\end{proof}
\begin{note}
Suppose that $(\sum_{i=2}^{n}k_i)+2=-kn$. Then the equalities 
\[ v_0+v_{n+1}+\sum_{i=2}^nk_iv_i=0,\ v_1+v_2+\cdots +v_n=0 \]
are equivalent to
\[ v_0+v_{n+1}+kv_1+\sum_{i=2}^n(k_i+k)v_i=0,\ v_1+v_2+\cdots +v_n=0, \]
and we have $k+\sum_{i=2}^{n}(k_i+k)=-2$. If we put $\bar{k}_i=k$ for 
$i=1$ and $\bar{k}_i=k_i+k$, then the projective space bundle 
$\Proj(\bigoplus_{i=2}^{n}\xi^{k_i})\oplus 1)$ is biholomorphic 
to $\Proj(\bigoplus_{i=1}^{n}\xi^{\bar{k}_i})$. Thus 
a complete non-singular toric variety $M$ of dimension $n$ with  
$c_1(M)$ divisible by $n$ can be written in the form
\[ M=\Proj(\bigoplus_{i=1}^{n}\xi^{\bar{k}_i})
 \ \text{with}\ \sum_i\bar{k}_i=-2. \]
The expression in \cite{Fuj} is given in this form.
\end{note}


\section{Appendix; orbifold elliptic genus}
In this section we recall the definitions of elliptic genus and of 
orbifold elliptic genus of an almost complex closed orbifold.
For the latter we follow \cite{BL2}. 
These genera can be defined also for stably almost complex orbifolds. 
For the sake of simplicity we shall confine ouselves to almost 
complex orbifolds. We also show that 
the genera $\var(\Delta(M),\V)$ and $\hatvar(\Delta(M),\V)$ of 
the pair of multi-fan $\Delta(M)$ and set of vectors $\V$ 
associated with an almost complex torus orbifold $M$ as was 
defined in Section 3 coincide with
those of $M$ in the sense of this section.

We first explain some facts about orbifolds needed to define
orbifold elliptic genus. We refer to \cite{HM} for basic
definitons and notations used here. 
Specifically $(V_x,U_x,H_x,p_x)$ will denote a reduced orbifold 
chart centered at $x\in M$, where $H_x$ is the isotropy group at 
$x$ and $p_x^{-1}(x)$ is a single point in $V_x$. We always assume
that the manifold $V_x$ is connected and small enough so that 
the fixed point set
of $h$, denoted by $V_x^h$, is also connected for each $h\in H_x$.
$p_x$ induces a homeomorphism $V_x/H_x\to U_x$. When 
$p_x$ has an obvious meaning in the context we shall omit $p_x$ 
and simply write $(V_x,U_x,H_x)$ for orbifold chart. 

Let $M$ be a connected closed orbifold. Kawasaki \cite{K1} defined 
an orbifold $\hatM$ and an orbifold map $\pi:\hatM\to M$ 
which plays an important
role in the index theorem of the Atiyah-Singer type for orbifolds. 
Some authors call each connected 
component of $\hatM$ by the name \emph{sector} of $M$. We shall call 
$\hatM$ the total sector 
of $M$. As a set $\hatM$ is defined by
\[ \hatM=\bigcup_{x\in M}Conj(H_x) \]
where $H_x$ is the isotorpy group at $x\in M$ and $Conj(H_x)$ denotes
the set of conjugacy classes of $H_x$.
The projectin $\pi$ is the obvious one. The orbifold structure of $\hatM$ is
given as follows. Let $(V_x,U_x,H_x,p_x)$ be a reduced orbifold chart 
of $M$ centered at $x$.  We put
\[ \hatV=\{(v,h)\in V_x\times H_x \mid hv=v \}. \] 
$H_x$ acts on $\hatV$ by
\[ h_1(v,h)=(h_1v,h_1hh_1^{-1}),\quad h_1\in H_x .\] 
Set $\hatU=\hatV/H_x$.
If $(v,h)\in \hatV$ and $y=p_x(v)\in U_x$, then the isotropy group $H_y$ 
is contained in $H_x$ and $h$ is contained in $H_y$
since we assumed that $V_x$ is small enough.
Therefore we get a map $\hatV\to \hatM$ sending $(v,h)$ to 
$[h]\in Conj(H_y)$, the conjugacy class of $h$ in $Conj(H_y)$. 
This map factors through $\hatU=\hatV/H_x\to \hatM$ which is injective.
This defines an orbifold chart of $\hatM$. 
Note that $\hatU$ is not always connected. In fact 
$\hatU=\bigcup_{[h]\in Conj(H_x)}V_x^h/C(h)$ is the decomposition into 
connected components, where $C(h)$ is the centralizer of 
$h$ in $H_x$. In particular $(V_x^h,V_x^h/C(h),C(h))$ is 
a reduced orbifold chart centered at $[h]\in Conj(H_x)\subset \hatM$. 
The projection $\pi:\hatM\to M$ is locally injective and is an orbifold map
as follows from the above definition. Each connected component $X$ of 
$\hatM$ is 
an orbifold. If $[h]$ is in the principal stratum of $X$, then the 
multiplicity $d(X)$ of $X$ is equal to the order $|C(h)|$ of $C(h)$.

The index theorem of Kawasaki \cite{K2} expresses the index of the Dirac 
operator twisted by a vector bundle as an integral of 
certain characteristic form over $\hatM$ where we assume $M$ to be
an almost complex orbifold. 
If a compact connected group $G$ acts on $M$ preserving the almost 
complex structure
then the index is defined as a virtual character of the group $G$,
and its value at $g\in G$ is expressed by Vergne's fixed 
point formula \cite{V}. Note that, if $G$ is connected and
$(V_x,U_x,H_x)$ is a reduced orbifold chart centered at $x\in M^G$ such that 
$U_x$ is $G$-invariant, then it can be shown that
some finite covering group $\tilde{G}$ of $G$ acts on $V_x$ in such a way that
the action commutes with the action of $H_x$ and covers the action of $G$ 
on $U_x$.  From this fact it follows that, 
if $G$ is a compact connected group, then 
the fixed point set $M^G$ is a suborbifold of $M$.

For the sake of simplicity we shall explain the fixed 
point formula only for torus actions. 
Let $M$ be an almost complex closed orbifold and $\xi$ an orbifold complex
vector bundle over $M$, both acted on by a torus $T$ topologically
generated by $g\in T$. The fixed point set $M^T$ is an almost complex
orbifold. Let $\hatM^T$ be the total sector of $M^T$ with 
projection $\pi:\hatM^T\to M^T$. We put $\hat{\xi}=\pi^*(\xi)$.
Let $(V_x,U_x,H_x)$ be a reduced orbifold chart of the 
orbifold $M$ centered at $x\in M^T$
having properties as explained above. We may suppose that the topological 
generator $g$ acts on $V_x$. Then it can be shown that its action commutes 
with that of $H_x$, and $(V_x^g,U_x^g,H_x)$ is a reduced orbifold
chart of the orbifold $M^T$ centered at $x$. We put
\[ \hatV^g =\{(v,h)\in V_x\times H_x \mid gv=v, hv=v\}.\] 
$H_x$ acts on $\hatV^g$ and $\hatU^g=\hatV^g/H_x$ maps homeomorphically 
onto an open subset of 
$\hatM^T$ as before. More precisely it can be shown that 
$\hatU^g=\bigcup_{[h]\in Conj(H_x)}(V_x^g)^h/C(h)$  
is the decomposition into connected components and 
$((V_x^g)^h,(V_x^g)^h/C(h),C(h))$ gives a reduced orbifold 
chart of $\hatM^T$ centered at $[h]$.

Let $F$ be a component of $M^T$ and $\hatF$ its total sector. 
We denote by $N_F$ the normal bundle of $F$ in $M$ and put 
$N_{\hatF}=\pi^*(N_F)$. Let $Y$ be a connected component 
of $\hatF$ and $N_Y$ the normal bundle of the immersion $\pi|Y:Y\to F$.
Suppose that $[h]\in Conj(H_x)$ lies in $Y$. 
The action of $g$ and $h$ on $V_x$ decomposes the 
tangent bundle $TV_x$ restricted to $(V_x^g)^h$ into 
a sum of eigen-bundles:
\[ TV_x|(V_x^g)^h=T(V_x^g)^h\oplus \bigoplus_\lambda W_\lambda\oplus 
\bigoplus_{\mu,\lambda}W_{\mu,\lambda}, \]
where $h$ acts on the fibers of $W_\lambda$ and $W_{\mu,\lambda}$ 
by multiplication by $\lambda$, and $g$ acts trivially on $W_\lambda$ 
and acts on the fibers of $W_{\mu,\lambda}$ by the multipliction by $\mu$. 
Since the action of $C(h)$ 
commutes with those of $g$ and $h$, we get orbifold vector bundles
$E_\lambda=W_\lambda/C(h)$ and $E_{\mu,\lambda}=W_{\mu,\lambda}/C(h)$
over $(V_x^g)^h/C(h)$. Thus, locally, we get the decomposition into a sum 
of eigen-bundles:
\begin{align*}
N_Y&=\bigoplus E_\lambda, \\
N_{\hatF}|Y&=\bigoplus E_{\mu,\lambda}.
\end{align*}
This decomposition is valid throughout the sector $Y$ since the 
eigenvalues $\lambda$ and $\mu$ are functions of conjugacy classes and it is 
constant on the connected space $Y$. 
We also have the decomposition $\hat{\xi}|Y=\bigoplus \xi_{\mu,\lambda}$
in a similar way.
For an orbifold complex vector bundle $E$ let $\Omega(E)$ denote the curvature
form of a connection on $E$ associated with invariant hermitian metrics on the 
tangent bundle of the base space and $E$. 
We put $\Gamma(E)=\frac{i}{2\pi}\Omega(E)$ and 
\begin{align*}
 \mathcal{T}(Y)&=\det\left(\frac{\Gamma(TY)}{1-e^{-\Gamma(TY)}}\right), \\
 \ch_g(\hat{\xi}|Y)&=\sum_{\mu,\lambda}\mu\lambda
 \tr e^{\Gamma(\xi_{\mu,\lambda})}, \\
 \mathcal{D}(N_Y)&=\prod_\lambda\det(1-\lambda^{-1}e^{-\Gamma(E_\lambda)}), \\
 \mathcal{D}(N_{\hatF}|Y)&=\prod_{\mu,\lambda}
 \det(1-(\mu\lambda)^{-1}e^{-\Gamma(E_{\mu,\lambda})}).
\end{align*}
Let $\partial$ denote the Dirac operator of $M$. Then Vergne's fixed point 
formula reads as follows:
\begin{equation}\label{eq:Vergne}
 \ind(\partial\otimes\xi)_g=\sum_F\sum_Y\frac1{d(Y)}\int_{Y}
 \frac{\mathcal{T}(Y)\ch_g(\hat{\xi}|Y)}
 {\mathcal{D}(N_Y)\mathcal{D}(N_{\hatF}|Y)},
\end{equation}
where the left hand side is the value of $\ind(\partial\otimes\xi)$ regarded
as a virtual character of $T$ at $g$, and the sums at the right hand side 
extend over all components $F$ of $M^T$ and all components $Y$ of $\hatF$. 
 
If $M$ is an almost complex orbifold, the 
elliptic genus $\var(M)$ of $M$ is defined as the index of the
Dirac operator $\partial$ twisted by the vector bundle 
\begin{equation}\label{eq:elliptic}
 \Lambda_{-\zeta}T^*M\otimes
\bigotimes_{k\geq 1}\left(\Lambda_{-\zeta q^k}T^*M
 \otimes\Lambda_{-\zeta^{-1}q^k}TM
    \otimes S_{q^k}T^*M\otimes S_{q^k}TM\right),
\end{equation}
where $TM$ and $T^*M$ are the tangent bundle and cotangent bundle of $M$
respectively, and $\Lambda_y$ and $S_y$ denote total exterior power and
symmetric power respectively. 

We next proceed to define the orbifold elliptic genus $\hatvar(M)$.
Let $X$ be a component of $\hatM$. It is also an almost complex
orbifold. Let $[h]\in Conj(H_x)\subset X$ be a point in $X$. 
Let $N_X$ denote the orbifold normal bundle of the 
immersion $\pi\vert X:X\to M$. 
As before we have the decomposition of $N_X$ into a sum of eigen-bundles
with respect to the action of $h$:
\[ N_X=\bigoplus_\lambda E_\lambda, \]
where $h$ acts on the fibers of $E_\lambda$ by multiplication by 
$\lambda=e^{2\pi\sqrt{-1}f_{X,\lambda}}\not=1$. We make the convention that
$0<f_{X,\lambda}<1$. We set 
$f_X=\sum_{\lambda}(\rank E_{\lambda})f_{X,\lambda}$, 
and
\[ \E_X=\bigotimes_{\lambda\not=1}\bigotimes_{k=1}^\infty\left(
 \Lambda_{-\zeta q^{f_{X,\lambda}}q^k}E_\lambda^*\otimes
 \Lambda_{-\zeta^{-1}q^{-f_{X,\lambda}}q^k}E_\lambda\otimes
 S_{q^{f_{X,\lambda}}q^k}E_\lambda^*\otimes 
 S_{q^{-f_{X,\lambda}}q^k}E_\lambda\right). \]

We now define
\begin{equation}\label{eq:orbelliptic}
 \begin{split}
 \hatvar(M)=& \zeta^{-\frac{n}{2}}\sum_{X\subset \hatM}\zeta^{f_X} \\
 &\ind\left(
 \partial\otimes\Lambda_{-\zeta}T^*X\otimes\bigotimes_{k=1}^\infty
 \left(\Lambda_{-\zeta q^k}T^*X\otimes\Lambda_{-\zeta^{-1}q^k}TX\otimes
 S_{q^k}T^*X\otimes S_{q^k}TX\right)\otimes\E\right). \end{split} 
\end{equation}
Note that, 
if we write $TX=E_1$ and $f_{X,1}=0$, then we can also write 
as follows:
\begin{equation}\label{eq:orbelliptic bis}
\begin{split}
 \hatvar(M)=&\zeta^{-\frac{n}{2}}\sum_{X\subset \hatM}\zeta^{f_X} 
 \ind(
 \partial\otimes\Lambda_{-\zeta}T^*X\otimes \\
 &\bigotimes_\lambda\bigotimes_{k=1}^\infty (
 \Lambda_{-\zeta q^{f_{X,\lambda}}q^k}E_\lambda^*\otimes
 \Lambda_{-\zeta^{-1}q^{-f_{X,\lambda}}q^k}E_\lambda\otimes
 S_{q^{f_{X,\lambda}}q^k}E_\lambda^*\otimes 
 S_{q^{-f_{X,\lambda}}q^k}E_\lambda ) ). \end{split} 
\end{equation}

Suppose that a torus $T$ acts on the orbifold $M$. Then 
\eqref{eq:elliptic} and \eqref{eq:orbelliptic} can be considered as defining
virtual characters of $T$ which are called equivariant 
elliptic genus and equivariant orbifold elliptic genus respectively.
We shall explicitly write down formulae for the equivariant elliptic 
genus $\var(M)$ and the 
equivariant orbifold elliptic genus $\hatvar(M)$, assuming that the fixed 
points are all isolated and the isotropy group at each fixed point is abelian
for the sake of simplicity. Write $M^T=\{P_j\}$. 
Let $(V_j,U_j,H_j,p_j)$ be a reduced orbifold chart centered at $P_j$, 
and let $\tilde{x}_j\in V_j$ be the unique point with 
$p_j(\tilde{x}_j)=P_j$. Note that the portion 
of $\hatM^T$ over $P_j$ is identified with the set $\{P_{j,h}\mid h\in H_j\}$ 
in this case, where $P_{j,h}=(\tilde{x}_j,h)\in \hat{V}_j$.

Let $g$ be a topological generator of $T$. The actions of 
$g$ and $H_j$ decompose
the tangent space $T_{\tilde{x}_j}V_j$ into a sum of 
one-dimensional eigen-spaces:
\begin{equation}\label{eq:eigen}
 T_{\tilde{x}_j}V_j=\bigoplus_i W_{i,j} 
\end{equation}
where $g$ acts on $W_{i,j}$ by multiplication by 
$e^{2\pi\img m_{i,j}z}$ and $h\in H_j$ 
acts by multiplication by $\chi_{i,j}(h)=e^{2\pi\img m_{i,j}^{H_j}}(h)$ where 
$\chi_{i,j}$ is a character of $H_j$.
We can also view \eqref{eq:eigen} as representing the decomposition of 
\[ N_{P_{j,h}}\oplus N_{\hatF}|P_{j,h}. \]
Then, under the above situation, 
we obtain the following formula from \eqref{eq:elliptic} and
\eqref{eq:Vergne}.
\begin{equation}\label{eq:elliptic formula}
 \var(M)_g=\sum_j\frac1{|H_j|}\sum_{h\in H_j}
             \prod_i\phi(-m_{i,j}z-m_{i,j}^{H_j}(h),\tau,\sigma). 
\end{equation}

Next let $X$ be a component of $\hatM$. We denote by $J_X$ the set 
of $j$ such that $P_j\in \pi(X)$. 
Then, for any $j\in J_X$, there is a unique $h\in H_j$ such that 
the point $P_{j,h}$ over $P_j$ 
lies in $X$. This $h$ will be denoted by $h_{X,j}$. 
We see that $W_{i,j}$ or 
rather $W_{i,j}/C(h)$ is tangent to $\pi(X)$ if and only 
if $\chi_{i,j}(h_{X,j})=1$. 
Moreover, if we write 
$\chi_{i,j}(h)=e^{2\pi\sqrt{-1}f_{i,j}(h)}$ with $0\leq f_{i,j}(h)<1$, 
then $f_{X}=\sum_if_{i,j}(h_{X,j})$. Note that this equality holds for 
any $j\in J_X$. The total sector
of the single point $P_{j,h}$ consisits of $\{P_{j,h'}\mid h'\in H_j\}$.
Then we obtain the following formula from
\eqref{eq:orbelliptic bis} and \eqref{eq:Vergne}.
\begin{equation}\label{eq:orbelliptic formula}
\hatvar(M)_g=\sum_{X\subset \hatM}\zeta^{f_X}\sum_{j\in J_X}
 \sum_{h'\in H_j}\frac1{|H_j|}
 \prod_i\phi(-m_{i,j}z+f_{i,j}(h_{X,j})\tau-m_{i,j}^{H_j}(h'),\tau,\sigma).
\end{equation}   

 For $P_j$ and $h\in H_j$ we set $f_{j,h}=\sum_if_{i,j}(h)$. 
The above formula \eqref{eq:orbelliptic formula} can be transformed in 
\begin{equation}\label{eq:orbelliptic formula2}
\hatvar(M)_g=\sum_j\sum_{(h_1,h_2)\in H_j\times H_j}
\frac{\zeta^{f_{j,h_1}}}{|H_j|}
 \prod_i\phi(-m_{i,j}z+f_{i,j}(h_1)\tau-\mu_{i,j}^{H_j}(h_2),\tau,\sigma).
\end{equation}   
In fact, if we set 
$ \mathcal{A}=\{(X,j)\mid X\ \text{is a component of $\hatM$},\ j\in J_X \} $
and $\mathcal{B}=\bigsqcup_jH_j$, then there is a bijection
$\rho:\mathcal{A}\to\mathcal{B}$ which sends 
$(X,j)\in\mathcal{A}$ to $h_{X,j}\in \mathcal{B}$, its
inverse being $H_j\ni h\mapsto (X,j)\in \mathcal{A}$ where $X$ is the 
component containing $P_{j,h}$. Moreover $f_X=f_{j,h}$ with 
this correspondence. 

Now let $M$ be an almost complex torus orbifold of dimension $2n$ 
acted on by a torus $T$.
Let $\Delta(M)=(\Sigma(M),C(M),w(M)^\pm)$ and 
$\V(M)=\{v_i\}_{i\in \Sigma(M)^{(1)}}$ be the multi-fan and the set 
of generating edge vectors associated with $M$. $\Delta(M)$ is 
a multi-fan in the lattice $L=\Hom(S^1,T)$.
We shall abbreviate the letter $M$ and simply write $\Delta$ for $\Delta(M)$ 
and so on. The fixed point set $M^T$ is the union of $M_I$ with $I\in \sigmn$. 
\begin{lemm}\label{lemm:isotropy}
If $P_j\in M^T$ lies in $M_I$, then the isotropy group $H_j$ of $P_j$ is 
isomorphic to $H_I=L/L_{I,\V}$.
\end{lemm}
\begin{proof}
Let $(V_j,U_j,H_j,p_j)$ be a reduced orbifold chart such that $U_j$ 
is invariant under the action of the torus $T$. The covering group 
$\tilde{T}_I=\prod_{i\in I}\tilde{S}_i$ of $T$ acts effectively
on $V_j$ where $\tilde{S}_i$ was introduced in Section 2. The 
isotropy group $H_j$ is isomorphic to the kernel $H$ of the natural
homomorphism $\tilde{T}_I\to T$ as was shown in \cite{HM}. 

On the other hand $L$ and $L_{I,\V}$ are identified with $\pi_1(T)$ 
and $\pi_1(\tilde{T}_I)$. Therefore the kernel $H$ is isomorphic to
$H_I=L/L_{I,\V}$. Hence $H_j$ is isomorphic to $H_I$.
\end{proof}
In the sequel we shall identify $H_j$ with $H_I$ where $P_j\in M_I$.
The decomposition
\eqref{eq:eigen} of the tangent space $T_{\tilde{x}_j}V_j$ into sum of
$1$-dimensional $\tilde{T}_I$ modules can be written 
in this case as
\begin{equation}\label{eq:eigen2}
 T_{\tilde{x}_j}V_j=\sum_{i\in I}t^{u_i^I}. 
\end{equation}
Hence $f_{i,j}(h_1)=f_{I,h_1,i}$ for $h_1\in H_j=H_I$ where 
$f_{I,h_1,i}$ is as in Section 3. It follows that 
$f_{j,h_1}=f_{I,h_1}$. Let $v(h_1)$ be the representative of 
$h_1\in H_I$ such that $\l u_i^I,v(h_1)\r =f_{I,h_1,i}$. 
Note also that $w^+(I)=\#\{P_j\in M_I\}$ and $w^-(I)=0$ in this case
since $M$ is an almost complex torus orbifold. 

From these observations and from \eqref{eq:elliptic formula} and 
\eqref{eq:orbelliptic formula2}
we obtain
 \begin{equation}\label{eq:multiformula2}
  \varv(M)=\sum_{I\in \sigmn}\frac{w(I)}{|H_I|}
 \sum_{h\in H_I}
 \prod_{i\in I}\phi(\l u_i^I, -zv-v(h)\r, \tau,\sigma), 
\end{equation}
and
\begin{equation}\label{eq:multiorbiformula3}
 \hatvarv(M)=\sum_{I\in \sigmn}\frac{w(I)}{|H_I|}
 \sum_{(h_1,h_2)\in H_I\times H_I} \zeta^{f_{I,h_1}}
 \prod_{i\in I}\phi(\l u_i^I, -zv+\tau v(h_1)-v(h_2)\r, \tau,\sigma).
\end{equation}

In order to get the formula for $\hatvarv(M)$ correspondig to 
\eqref{eq:orbelliptic formula}, we look for the pair $(X,j)\in \mathcal{A}$ 
which corresponds via $\rho$ to $P_{j,h}$ with $P_j\in M_I$ and $h\in H_I$ .
As was remarked just before \eqref{eq:orbelliptic formula}, $t^{u_i^I}$ 
in \eqref{eq:eigen2} is tangent to $X$ if and only if $f_{I,h,i}=
f_{i,j}(h)=0$. This means that $\pi(X)$ is contained in $M_K$ and 
is one of its component, where 
$K=\{i\in I\mid f_{I,h,i}\not=0\}$. Moreover $h$ is contained in $\hatH_K$,
where $\hatH_K$ is defined in Section 3 and is characterized by \eqref{eqn:0}.
But \eqref{eqn:0} is equivalent to
\[ \hatH_K=\{h\in H_I\mid f_{I,h,i}=0\ 
\text{for}\ i\in I\setminus K\ \text{and}\ f_{I,h,i}\not=0\ 
\text{for}\ i\in K \}.\]
In particular $f_{K,h}=f_{I,h}$. Thus \eqref{eq:multiorbiformula3} can be 
rewritten in the following form.
\begin{equation}\label{eq:multiorbiformula4}
\begin{split}
 \hatvarv(M) 
 = &\sum_{k=0}^n\sum_{K\in\sigmk,h_1\in \hat{H}_K}
  \zeta^{f_{K,h_1}} 
  \sum_{I\in \Sigma_K^{(n-k)}}\frac{w(I)}{|H_I|}\cdot \\
   &\sum_{h_2\in H_I}\left(
  \prod_{i\in I\setminus K}\phi(\l u_i^I,-zv-v(h_2)\r,\tau,\sigma)
  \prod_{i\in K}\phi(\l u_i^I,-zv+\tau v(h_1)-v(h_2)\r,\tau,\sigma) 
  \right).
\end{split}
\end{equation}
  
\begin{rema}
So far $M$ is assumed to be an almost complex torus orbifold. Even 
if $M$ is a stably almost complex torus orbifold the formulae 
\eqref{eq:multiformula2}, \eqref{eq:multiorbiformula3} 
and \eqref{eq:multiorbiformula4} are valid. We have only to count the points
$P_j\in M_I$ with sign and the resulting multiplicity is $w(I)$. 
The formulae \eqref{eq:varv} for $\varv(\Delta,\V)$ and 
\eqref{eq:hatvarv} and  \eqref{eq:hatvarv2} 
for $\hatvarv(\Delta,\V)$ given in Section 3 
are modelled on \eqref{eq:multiformula2}, \eqref{eq:multiorbiformula3} and 
\eqref{eq:multiorbiformula4} respectively.
\end{rema}


\providecommand{\bysame}{\leavevmode\hbox to3em{\hrulefill}\thinspace}

\end{document}